\documentclass{amsart}
\usepackage{amsmath,amssymb,amscd,latexsym}
\usepackage[all]{xy}

\newcommand{\C}{{\mathbb{C}}}

\newcommand{\N}{{\mathbb{N}}}

\newcommand{\R}{{\mathbb{R}}}

\newcommand{\Z}{{\mathbb{Z}}}

\newcommand{\Ah}{{\mathcal A}}

\newcommand{\Ch}{{\mathcal C}}
\newcommand{\Dh}{{\mathcal D}}

\newcommand{\Fh}{{\mathcal F}}
\newcommand{\Gh}{{\mathcal G}}

\newcommand{\Kh}{{\mathcal K}}

\newcommand{\Oh}{{\mathcal O}}

\newcommand{\Qh}{{\mathcal Q}}

\newcommand{\Uh}{{\mathcal U}}

\newcommand{\Zh}{{\mathcal Z}}

\newcommand{\ad}{\mathrm{ad}\,}

\newcommand{\au}{\approx_{\mathrm{au}}}
\newcommand{\asu}{\approx_{\mathrm{asu}}}
\newcommand{\sau}{\approx_{\mathrm{sau}}}
\newcommand{\sasu}{\approx_{\mathrm{sasu}}}

\newcommand{\be}{\mathbf{1}}

\newcommand{\dist}{\mathrm{dist}}

\newcommand{\ev}{\mathrm{ev }}

\newcommand{\id}{\mathrm{id}}
\newcommand{\Inv}{\mathrm{Inv}}

\newcommand{\Prim}{\mathrm{Prim}}

\newcommand{\ue}{\approx_{\mathrm{u}}}

\newcounter{number}[section]

\newenvironment{nummer}{\refstepcounter{number}{\noindent\arabic{section}.\arabic{number}}}{}

\newcommand{\bn}{\noindent \begin{nummer} \rm}
\newcommand{\en}{\end{nummer}}

\newenvironment{ntheorem}{\noindent {\sc Theorem:} \it}{}
\newenvironment{nlemma}{\noindent {\sc Lemma:} \it}{}
\newenvironment{nprop}{\noindent {\sc Proposition:} \it}{}
\newenvironment{ndefn}{\noindent {\sc Definition:} \it}{}
\newenvironment{ncor}{\noindent {\sc Corollary:} \it}{}

\newenvironment{nremark}{\noindent {\sc Remark:}}{}
\newenvironment{nremarks}{\noindent {\sc Remarks:}}{}
\newenvironment{nexamples}{\noindent {\sc Examples:} }{}

\newenvironment{nnotation}{\noindent {\sc Notation:} }{}

\newenvironment{nproof}{\noindent {\sc Proof:}}{\mbox{}\hfill 
\rule[-.2ex]{.25em}{1.8ex}}

\begin{document}

\title[Localizing the Elliott conjecture]{{\sc Localizing the Elliott conjecture\\
at strongly self-absorbing $C^{*}$-algebras\\
\footnotesize -- with an appendix by Huaxin Lin --}}

\author{Wilhelm Winter}
\address{Mathematisches Institut der Universit\"at M\"unster\\
Einsteinstr.\ 62\\ D-48149 M\"unster}

\email{wwinter@math.uni-muenster.de}

\date{\today}
\thanks{{\it Supported by:}  Deutsche 
Forschungsgemeinschaft (SFB 478)}

\setcounter{section}{-1}

\begin{abstract}
We formally introduce the concept of localizing the Elliott conjecture at a given strongly self-absorbing $C^{*}$-algebra $\Dh$; we also explain how the known classification theorems for nuclear $C^{*}$-algebras fit into this concept. As a new result in this direction, we employ recent results of Lin to show that (under a mild $K$-theoretic condition) the class of separable, unital, simple $C^{*}$-algebras with locally finite decomposition rank and UCT, and  for which projections separate traces, satisfies the Elliott conjecture localized at the Jiang--Su algebra $\mathcal{Z}$.

Our main result is formulated in a  more general way; this allows us to outline a strategy to possibly remove the trace space condition as well as the $K$-theory restriction entirely. When regarding both our result and the recent classification theorem of Elliott, Gong and Li as generalizations of the real rank zero case, the two approaches are perpendicular in a certain sense. The strategy to attack the general case aims at combining these two approaches.

Our classification theorem covers simple ASH algebras for which projections separate traces (and the $K$-groups of  which have finitely generated torsion part); it does, however, not at all  depend on an inductive limit structure. Also, in the monotracial case it does not rely on the existence or absence of projections in any way. In fact, it is the first such result which, in a natural way, covers all known unital, separable, simple, nuclear and stably finite $C^{*}$-algebras of real rank zero (the $K$-groups of  which have finitely generated torsion part) as well as the (projectionless) Jiang--Su algebra itself. 
\end{abstract}

\maketitle

\section{Introduction}

It is the aim of the Elliott program to classify separable nuclear $C^{*}$-algebras by their $K$-theory data, see  \cite{Ell:classprob} and \cite{Ror:encyc} for an overview. While there are some striking results in the non-simple case (\cite{KirRor:pi}, \cite{KirRor:pi2}, and in particular \cite{Kir:Michael} and \cite{Kir:fields}), the Elliott  conjecture has most successfully been verified for classes of simple $C^{*}$-algebras (cf.\ \cite{Ell:simp}, \cite{EllGongLi:simple_AH}, \cite{KirPhi:classI}, \cite{KirPhi:classII}) -- and this is the case we shall mostly be dealing with in the present notes. 

In recent years, the notion of strongly self-absorbing $C^{*}$-algebras has become of particular importance for the Elliott program. A $C^{*}$-algebra $\Dh \neq \C$ is strongly self-absorbing if $\Dh$ is isomorphic to $\Dh \otimes \Dh$ such that the isomorphism is approximately unitarily equivalent to the first factor embedding, cf.\ \cite{TomsWinter:ssa}. A $C^{*}$-algebra $A$ is $\Dh$-stable if $A \cong A \otimes \Dh$. The list of known strongly self-absorbing $C^{*}$-algebras is relatively short (cf.\ \ref{D-examples} below); it consists of the Cuntz algebras $\Oh_{2}$, $\Oh_{\infty}$, UHF algebras of infinite type,  tensor products of $\Oh_{\infty}$ with UHF algebras of infinite type and the Jiang--Su algebra $\Zh$. Kirchberg's celebrated classification results may clearly be regarded as classification up to $\Oh_{\infty}$-stability or up to $\Oh_{2}$-stability, respectively. It was shown in \cite{TomsWinter:Zash} that essentially all currently known classification theorems in the framework of the Elliott program may  be interpreted  as classification up to $\Zh$-stability. There are also some classification type results up to UHF-stability, cf.\ \cite{Dad:UHFhomotopy}. In the present article, we formalize the concept of classification up to $\Dh$-stability  by introducing the notion of localizing the Elliott conjecture  at a given strongly self-absorbing $C^{*}$-algebra $\Dh$.   The counterexamples to the Elliott conjecture provided by R{\o}rdam and by Toms (based on Villadsen's ideas, cf.\ \cite{Ror:simple}, \cite{Toms:classproblem}, \cite{Vil:perforation} and \cite{Vil:sr=n}) even suggest that $\Dh$-stability for at least some strongly self-absorbing $\Dh$ will be necessary for classification.

The algebras $\Oh_{\infty}$ and $\Zh$ are of particular relevance, since they are both $KK$-equivalent to the complex numbers and therefore $\Oh_{\infty}$-stability and $\Zh$-stability are the least restrictive on $K$-theory. In fact, we think of $\Zh$ as stably finite analogue of $\Oh_{\infty}$; and given that Kirchberg--Phillips classification of purely infinite $C^{*}$-algebras really covers simple, nuclear and $\Oh_{\infty}$-stable $C^{*}$-algebras, it is natural to ask for similar statements in the purely finite (cf.\ \cite[3.3]{TomsWinter:VI}) case, i.e., for simple, nuclear,  finite and  $\Zh$-stable $C^{*}$-algebras.  The results of \cite{Winter:Z-class} and \cite{Winter:lfdr} (based on ideas from \cite{KirWinter:dr} and  \cite{Winter:fintopdim}) are of this nature. They cover $\Zh$-stable simple $C^{*}$-algebras of (locally) finite decomposition rank with real rank zero. One should note that locally finite decomposition rank is a fairly weak structural hypothesis on stably finite $C^{*}$-algebras (it holds for all ASH, i.e., approximately subhomogeneous algebras and  in particular for the pathological examples of Toms and Villadsen). Real rank zero (cf.\ \cite{BroPed:realrank}), however, is a serious restriction, although it is satisfied by many interesting and important examples, such as irrational rotation algebras.  It implies the existence of an abundance of projections;   it is  therefore no surprise that in the real rank zero case it is generally easier to extract information about the algebra from its ($K$-theory based) Elliott invariant.  Also, many of the earlier results in the subject deal exclusively with the real rank zero case, cf.\ \cite{Ell:rrzeroI}, \cite{EllGong:rrzeroII}, \cite{DadGong:class-rr0}.

It is of course natural to try to generalize these  results to a non real rank zero situation. For quite some time, Elliott--Gong--Li classification of simple AH algebras (with very slow dimension growth, cf.\ \cite{EllGongLi:simple_AH})  remained the only significant result for a natural and substantial class of simple $C^{*}$-algebras not necessarily of real rank zero (despite a number of very interesting yet somewhat scattered results which do not require the existence of projections, but only cover inductive limits of special building blocks, cf.\ \cite{Mygind:torsionclass} and \cite{JiaSu:Z}). The drawback in this otherwise major leap forward was, that it only covers algebras which a posteriori turned out to be approximately divisible (\cite{EllGongLi:apprdiv}), hence still contain relatively many projections. By \cite{BlaKumRor:apprdiv}, approximately divisible $C^{*}$-algebras for which projections separate traces have real rank zero. In our Corollary~\ref{fewtraces-lfdr-cor} below, we shall ask for this second condition rather than approximate divisibility to generalize the real rank zero case. In this sense, our approach is perpendicular (at least outside the class of real rank zero algebras) to that of Elliott, Gong and Li. 

Our classification result covers all known separable, unital, simple and nuclear $C^{*}$-algebras with finitely generated $K$-groups which are finite and have real rank zero; it also covers the Jiang--Su algebra itself, as well as certain projectionless crossed products considered by Connes (see \ref{classified-examples} below). 

Corollary~\ref{fewtraces-lfdr-cor} also has some other interesting consequences. For instance, using \cite{LinPhi:mindifflimits} we can make  progress on the classification problem for uniquely ergodic minimal dynamical systems (Corollary~\ref{uniquely-ergodic-cor}). Moreover, using recent work of Dadarlat and Toms, we  characterize the Jiang--Su algebra as the uniquely determined strongly self-absorbing and projectionless ASH algebra (Corollary~\ref{Z-ASH-characterization}).

We have split the way towards Corollary~\ref{fewtraces-lfdr-cor} in two parts: The first one is Theorem~\ref{main-result}, which shows how existing classification results for $C^{*}$-algebras which absorb UHF algebras can be used to derive classification for $\Zh$-stable $C^{*}$-algebras. The second part is the actual proof of Corollary~\ref{fewtraces-lfdr-cor}; it merely consists of verifying the (rather abstract) hypotheses of \ref{main-result}. This, however, will follow from earlier results by H.\ Lin and by the author, and from the appendix to this paper written by H.\ Lin. At this point, it also becomes clear that the only purpose of the additional conditions in \ref{fewtraces-lfdr-cor} is to make these real rank zero results available.

The next problem will be to try to  remove the trace space condition as well as the $K$-theory restriction of \ref{fewtraces-lfdr-cor} entirely. We believe that our approach (i.e., Theorem~\ref{main-result} together with the proof of Corollary~\ref{fewtraces-lfdr-cor}) paves the road for such an attempt; in fact, we outline a  strategy in \ref{outlook} below. Essentially, it suggests to use \ref{main-result} as a black box and verify its hypotheses for a larger class than that of \ref{fewtraces-lfdr-cor}. It is particularly tempting that this strategy uses existing classification theorems (with many projections) rather than just the (sometimes extremely complicated) methods developed to handle  these cases.

The paper is organized as follows: In Section 1, we recall some versions of approximate unitary equivalence as well as the notion of and some facts about strongly self-absorbing $C^{*}$-algebras. In Section 2 we introduce the concept of localizing the Elliott conjecture at a strongly self-absorbing $C^{*}$-algebra.   In Section 3, we study generalized dimension drop intervals and inductive limits of such.  The notion of unitarily suspended isomorphisms between tensor products with generalized dimension drop intervals is analyzed in Section 4. Here, we also provide  a technical key step for the proof of \ref{main-result}, namely, we construct an asymptotic unitary intertwining between $\Zh$-stable $C^{*}$-algebras from a unitarily suspended isomorphism.  In Section 5 we collect some observations on the $K$-theory of tensor products with UHF algebras.  In Section 6 we recall Lin's notion of tracial rank zero $C^{*}$-algebras as well as a result from \cite{Lin:existence}.  Our main result about localizing the Elliott conjecture at the Jiang--Su algebra  is Theorem~\ref{main-result};  we  derive its corollaries and applications in Section 8. Finally, we outline a strategy of how to possibly remove the trace space condition and $K$-theory restriction from Corollary~\ref{fewtraces-lfdr-cor} in Section 9.

Part of this work was carried out while the author was visiting Copenhagen University and the Centre de Recerca Matem{\`a}tica in Barcelona. I am indebted to both institutions, and in particular to S{\o}ren Eilers and Francesc Perera, for their kind hospitality. I am also grateful to Mikael R{\o}rdam, as well as to Marius Dadarlat, for inspiring conversations on the classification program in general, and  on the potential use of $\Zh$-stability and UHF-absorbtion in particular.  Finally, I would like to thank Huaxin Lin for pointing out a mistake in an earlier version of the proof of Corollary~\ref{fewtraces-lfdr-cor}, and for providing a way to fill that gap.

\section{Strongly self-absorbing $C^{*}$-algebras}

Let us recall some versions of approximate unitary equivalence as well as the notion of strongly self-absorbing $C^{*}$-algebras.

\bn
\begin{ndefn}
Let $A$, $B$ be separable  $C^{*}$-algebras and $\varphi,\psi: A \to B$ be $*$-homomorphisms. Suppose that $B$ is unital.
\begin{enumerate}
\item[(i)] We say $\varphi$ and $\psi$ are approximately unitarily equivalent, $\varphi \au \psi$, if there is a sequence $(u_{n})_{n \in \N}$ of unitaries in  $B$ such that
\[
\|u_{n} \varphi(a) u_{n}^{*} - \psi(a)\| \stackrel{n \to \infty}{\longrightarrow} 0
\]
for every $a \in A$.
If all $u_n$ can be chosen  to be in $\Uh_{0}(B)$, the connected
component of $\be_B$ of the unitary group $\Uh(B)$, then we say that
$\varphi$ and $\psi$ are strongly approximately unitarily equivalent,
written $\varphi \sau \psi$.
\item[(ii)] We say $\varphi$ and $\psi$ are asymptotically unitarily equivalent, $\varphi \asu \psi$, if there is a norm-continuous path $(u_{t})_{t \in [0,\infty)}$ of unitaries in  $B$ such that
\[
\|u_{t} \varphi(a) u_{t}^{*} - \psi(a)\| \stackrel{t \to \infty}{\longrightarrow} 0
\]
for every $a \in A$. If one can arrange that $u_0=\be_B$ and hence ($u_t\in\Uh_{0}(B)$ for all $t$), then we say that
$\varphi$ and $\psi$ are strongly asymptotically unitarily equivalent,
written $\varphi \sasu \psi$.
\end{enumerate}
\end{ndefn}
\en

\bn
The concept of strongly self-absorbing $C^{*}$-algebras was formally introduced in \cite[Definition~1.3]{TomsWinter:ssa}:
 
\begin{ndefn}
\label{def:ssa}
(i) A separable unital $C^{*}$-algebra $\Dh$ is strongly self-absorbing, if $\Dh \neq \C$ and there is an isomorphism $\varphi: \Dh \to \Dh \otimes \Dh$ such that $\varphi \au \id_{\Dh}\otimes \be_{\Dh}$.

(ii) A $C^{*}$-algebra $A$ is said to be $\Dh$-stable, if $A \cong A \otimes \Dh$.
\end{ndefn}
\en

\bn
\begin{nremark}
Strongly self-absorbing $C^{*}$-algebras are known to be simple and nuclear; moreover, they are either purely infinite or stably finite. $\Dh$-stability has good permanence properties, cf.\ \cite{TomsWinter:ssa}, \cite{Kir:spi}, \cite{Kir:CentralSequences}, \cite{HirshbergRordamWinter:D-stable}, \cite{HirshbergWinter:rokhlin-ssa}, \cite{HirshbergWinter:D-permutations} and \cite{DadWinter:KK-D-topology}.
\end{nremark}
\en

\bn
\label{D-examples}
\begin{nexamples}
The only known examples of strongly self-absorbing $C^{*}$-algebras are UHF algebras (\cite{Gli:UHF}) of infinite type (i.e., every prime number that occurs in the respective supernatural number occurs with infinite multiplicity -- in this case, we also say the supernatural number is of infinite type), the Cuntz algebras $\Oh_{2}$ and $\Oh_{\infty}$ (\cite{Cuntz:On}), the Jiang--Su algebra $\Zh$ (\cite{JiaSu:Z}) and tensor products of $\Oh_{\infty}$ with UHF algebras of infinite type. All these examples are known to be $K_{1}$-injective, i.e., the natural map $\Uh(\Dh)/\Uh_{0}(\Dh) \to K_{1}(\Dh)$ is injective.  Among UHF algebras of infinite type, the universal UHF algebra $\Qh$ will be of particular interest for us.
\end{nexamples}
\en

\section{Localizing Elliott's conjecture}

Below we recall the Elliott conjecture and formally introduce the concept of localizing it at a given strongly self-absorbing $C^{*}$-algebra $\Dh$. We also recall some existing classification results which support this point of view.

\bn
Roughly speaking, the Elliott conjecture says that separable, nuclear $C^{*}$-algebras are  classified by $K$-theoretic data. The proposed classifying invariant tends to be fairly complicated in general; however, when restricting to specific subclasses of nuclear $C^{*}$-algebras, it becomes significantly easier. In these notes, we are mainly interested in the simple, unital and stably finite case.

\begin{ndefn} 
\cite[Section~2.2]{Ror:encyc} For a separable, simple, unital and nuclear $C^{*}$-algebra $A$ we define its Elliott invariant to be
\[
\Inv A:= (K_{0}(A),K_{0}(A)_{+},[\be_{A}], K_{1}(A), T(A), r_{A}: T(A) \to S(K_{0}(A))).
\]
\end{ndefn}
Here, $ (K_{0}(A),K_{0}(A)_{+},[\be_{A}])$ is the preordered $K_{0}$-group, $K_{1}(A)$ the $K_{1}$-group, $T(A)$ the Choquet simplex of tracial states, $S(K_{0}(A))$ the Choquet simplex of states (i.e., positive unit-preserving group homomorphisms from $ (K_{0}(A),K_{0}(A)_{+},[\be_{A}])$ to $(\R,\R_{+},1)$) and  $r_{A}: T(A) \to S(K_{0}(A))$ the canonical affine map given by $r_{A}(\tau)([p])=\tau(p)$. 
\en

\bn
\begin{nremarks}
(i) The invariant as defined above is best suited to deal with stably finite $C^{*}$-algebras. In this case, $(K_{0}(A),K_{0}(A)_{+},[\be_{A}])$ in fact is an  ordered group with distinguished order unit in the sense of \cite[Definition~1.1.14]{Ror:encyc}. 

(ii) For purely infinite $C^{*}$-algebras, we have  $K_{0}(A)= K_{0}(A)_{+}$ (\cite[Proposition~2.2.2]{Ror:encyc}), and $T(A)=\emptyset$. So in this case, $\Inv A$ reduces to  $(K_{0}(A),[\be_{A}],K_{1}(A))$. 
\end{nremarks}
\en

\bn
\label{morphisms-of-invariants}
\begin{ndefn}
If $A$ and $B$ are separable, simple, unital and nuclear $C^{*}$-algebras, then by a morphism
\[
\Lambda: \Inv A \to \Inv B
\]
of invariants we mean a pair $\Lambda=(\lambda,\Delta)$, where 
\[
\lambda: (K_{0}(A),K_{0}(A)_{+},[\be_{A}], K_{1}(A)) \to (K_{0}(B),K_{0}(B)_{+},[\be_{B}], K_{1}(B))
\]
is a morphism of (pre-)ordered abelian groups and
\[
\Delta: T(B) \to T(A)
\]
is an affine map, such that the diagram
\begin{equation}
\label{24}
\xymatrix{
T(A) \ar[r]^{r_{A}} & S(K_{0}(A))  \\
T(B) \ar[r]^{r_{B}} \ar[u]^{\Delta} & S(K_{0}(B)) \ar[u]_{\lambda^{S}}
}
\end{equation}
(where $\lambda^{S}$ is the obvious map induced by $\lambda$) commutes.
\end{ndefn}
\en

\bn
\label{d-EC}
\begin{ndefn}
We say a class $\Ah$ of separable, simple, unital and nuclear $C^{*}$-algebras satisfies the Elliott conjecture (or  (EC), for short), if the following holds:

If $A$ and $B$ are $C^{*}$-algebras in $\Ah$ and if 
\[
\Lambda: \Inv A \to \Inv B
\]
is an isomorphism of invariants, then there is an isomorphism 
\[
\varphi: A \to B
\]
such that 
\[
\Inv(\varphi) = \Lambda.
\]
\end{ndefn}

The original form of the Elliott conjecture (in the unital case) states that  the class of \emph{all}  separable, simple, unital and nuclear $C^{*}$-algebras satisfies (EC), cf.\ \cite[Conjecture~2.2.5]{Ror:encyc}. 
\en

\bn
\label{d-uniqueness}
In many situations, it is important to be able to compare maps which agree on the invariant; in particular, one wants to know whether two such maps are approximately unitarily equivalent. We formalize this property as follows:

\begin{ndefn}
We say a class $\Ah$ of separable, simple, unital and nuclear $C^{*}$-algebras satisfies  $\Inv$-uniqueness, if the following holds:

If $A$ and $B$ are $C^{*}$-algebras in $\Ah$ and if
\[
\varphi_{1},\varphi_{2}: A \to B
\]
are $*$-isomorphisms such that
\[
\Inv(\varphi_{1}) = \Inv(\varphi_{2}),
\]
then 
\[
\varphi_{1} \au \varphi_{2}.
\]
If $\Ah$ consists of algebras in the UCT class, we say $\Ah$ satisfies $KL$-uniqueness, if in the above situation $KL(\varphi_{1})=KL(\varphi_{2})$ implies $\varphi_{1} \au \varphi_{2}$. Here, $KL$ denotes the bivariant functor introduced by R{\o}rdam in \cite{Ror:infsimple}, cf.\ also \cite[2.4.8]{Ror:encyc}.
\end{ndefn}
\en

\bn
\label{d-localized-EC}
\begin{ndefn}
Let $\Dh$ be a strongly self-absorbing $C^{*}$-algebra. We say a class $\Ah$ of separable, simple, unital and nuclear $C^{*}$-algebras satisfies the  Elliott conjecture \emph{localized at $\Dh$}, if the class
\[
\Ah^{\Dh}:= \{A \otimes \Dh \mid A \in \Ah\} 
\]
satisfies (EC) in the sense of \ref{d-EC}.

Similarly, we say $\Ah$ satisfies $\Inv$-uniqueness (or $KL$-uniqueness, respectively) localized at $\Dh$, if $\Ah^{\Dh}$ satisfies $\Inv$-uniqueness (or $KL$-uniqueness, respectively).
\end{ndefn}
\en

\bn 
The perhaps most spectacular manifestation of the Elliott conjecture is the classification of Kirchberg algebras  obtained by Kirchberg and by Phillips; see \cite{Kir:fields} and \cite{Phi:class}. An important ingredient of the proof is Kirchberg's first tensor product theorem, which states that a simple, separable and nuclear $C^{*}$-algebra is purely infinite if and only if it is $\Oh_{\infty}$-stable.  Kirchberg--Phillips classification may therefore clearly be interpreted as classification up to $\Oh_{\infty}$-stability. We state here only the  case of unital Kirchberg algebras satisfying the universal coefficients theorem (UCT).

\begin{ntheorem} 
Let $\Ah$ be the class of separable, simple, unital and nuclear $C^{*}$-algebras satisfying the UCT. Then, $\Ah$ satisfies the Elliott conjecture and $\Inv$-uniqueness, both localized at $\Oh_{\infty}$.
\end{ntheorem}
\en

\bn
In the stably finite case, there is a number of  classification results which so far are nonetheless somewhat incomplete. Below  we combine \cite[Theorem~2.1]{Winter:lfdr} with   \cite[Theorem~5.1]{Lin:TAFduke} (for the existence part) and with \cite[Theorem~1.1]{Dad:simpleTAFmorphisms} (for the uniqueness part). It is fairly general in the sense that it does not rely on an inductive limit structure in any way (locally finite decomposition rank is a \emph{local} condition); however, it only deals with the real rank zero case, i.e., it requires the existence of many projections.

\begin{ntheorem}
Let $\Ah$ be the class of separable, simple, unital $C^{*}$-algebras with locally finite decomposition rank, real rank zero, and satisfying the UCT.  Then, $\Ah$ satisfies the   Elliott conjecture and $KL$-uniqueness, both  localized at $\Zh$.
\end{ntheorem}
\en

\bn
We also wish to mention a result of M.~Dadarlat and the author which may be interpreted as a classification result  up to $\Dh$-stability. (However, since  it deals with non-simple $C^{*}$-algebras, it formally does not quite fit into the framework of Definition~\ref{d-localized-EC}.)

\begin{ntheorem} \cite[Theorem~0.1]{DadWinter:trivial-fields}
Let $A$ be a separable, unital $\Ch(X)$-algebra over a finite-dimensional compact metrizable space $X$. Suppose that all the fibres of $A$ are isomorphic to the same strongly self-absorbing $K_{1}$-injective $C^{*}$-algebra $\Dh$. Then, $A$ and $\Ch(X) \otimes \Dh$ are isomorphic as $\Ch(X)$-algebras.
\end{ntheorem}
\en

\section{Generalized dimension drop intervals}

In this section we study generalized dimension drop intervals and ways of writing the Jiang--Su algebra as limits of such.

\bn
\begin{ndefn}
Let $\mathfrak{p}$ and $\mathfrak{q}$ be supernatural numbers. We define the generalized dimension drop interval $Z_{\mathfrak{p},\mathfrak{q}}$ to be the $C^{*}$-algebra
\[
Z_{\mathfrak{p},\mathfrak{q}}:= \{f \in \Ch([0,1], M_{\mathfrak{p}} \otimes M_{\mathfrak{q}}) \mid f(0) \in M_{\mathfrak{p}} \otimes \be_{M_{\mathfrak{q}}}, \, f(1) \in \be_{M_{\mathfrak{p}}} \otimes M_{\mathfrak{q}} \}. 
\]
\end{ndefn}
We shall regard $Z_{\mathfrak{p},\mathfrak{q}}$ (and any tensor product with it) as $\Ch([0,1])$-algebra with the obvious central embedding of $\Ch([0,1])$. 
\en

\bn
It follows from \cite[Section~3]{JiaSu:Z} that prime dimension drop intervals are $KK$-equivalent to $\C$.  Writing a generalized dimension drop interval as an inductive limit of ordinary dimension drop intervals in the obvious way, it is a straightforward observation that the connecting maps are compatible with the $KK$-equivalences to $\C$.  It is then clear that the statement carries over to generalized dimension drop intervals, i.e., we have: 

\label{KK-C-Z-p-q}
\begin{nprop}
Let $\mathfrak{p}$ and $\mathfrak{q}$ be supernatural numbers which are relatively prime.  Then, $Z_{\mathfrak{p},\mathfrak{q}}$ is $KK$-equivalent to $\C$.
\end{nprop}
\en

\bn
\label{p-dimension-drop-inclusions-au}
\begin{nprop}
Let $p$, $q$ be natural numbers. 

(i) Suppose $p' \in \N$ divides $p$ and $q' \in \N$ divides $q$. Let $\beta, \gamma: Z_{p',q'} \to Z_{p,q}$ be unital  $\Ch([0,1])$-maps. Then, $\beta \au \gamma$.

(ii) Set $D  :=  \{f \in \Ch([0,1]^{2},M_{p} \otimes M_{q}) \mid  f(0,0) \in M_{p} \otimes \be_{q}, \, f(1,1) \in \be_{p} \otimes M_{q}\}$ and define embeddings
\[
\zeta_{x},\zeta_{y}: Z_{p,q} \to D 
\]
by
\begin{equation}
\label{8bb}
\zeta_{x}(f)(s,t) := f(s) \mbox{ and } \zeta_{y}(f)(s,t):= f(t). 
\end{equation}
There is a $\Ch([0,1]^{2})$ map 
\[
\alpha: D \to Z_{p,q} \otimes Z_{p,q}
\]
and for any such map we have
\begin{equation}
\label{1bb}
\alpha \circ \zeta_{x} \au \id_{Z_{p,q}} \otimes \be_{Z_{p,q}}
 \mbox{ and } 
\alpha \circ \zeta_{y} \au \be_{Z_{p,q}} \otimes \id_{Z_{p,q}}.
\end{equation}
\end{nprop}

\begin{nproof}
(i) Let $\iota_{p,q}: Z_{p,q} \to \Ch([0,1], M_{p} \otimes M_{q})$ denote the canonical embedding. Using that $\beta$ and $\gamma$ are both unital, it is not too hard to show that $\iota_{p,q} \circ \beta \au \iota_{p,q} \circ \gamma$ via  unitaries $u^{(n)} \in \Ch([0,1], M_{p} \otimes M_{q})$. We may even assume that each $u_{n}$ is constant on small intervals $[0,\varepsilon_{n}]$ and $[1-\varepsilon_{n},1]$, and that $u_{i}^{(n)}$ \emph{exactly} intertwines  $(\iota_{p,q} \circ \beta)_{i}$ and $( \iota_{p,q} \circ \gamma)_{i}$ for $i=0,1$.  

Similarly, we have $\beta_{0} \ue \gamma_{0}$, and $\beta_{1} \ue \gamma_{1}$, respectively, i.e., they are exactly unitarily equivalent.  Let $v_{0} \in M_{p} \otimes \be_{q}$ and $v_{1} \in \be_{p} \otimes M_{q}$ be unitaries implementing these  equivalences.  But then, $v_{0}^{*}u^{(n)}_{0} \in (\beta_{0}(M_{p'} \otimes \be_{q'}))' \subset M_{p} \otimes M_{q}$, whence there is a path $(s^{(n)}_{t})_{t \in [0,\varepsilon]}$ of unitaries in  the commutant $(\beta_{0}(M_{p'} \otimes \be_{q'}))' \subset M_{p} \otimes M_{q}$ such that $s^{(n)}_{\varepsilon}=v_{0}^{*} u^{(n)}_{0}$ and $s^{(n)}_{0} = \be_{p} \otimes \be_{q}$. A similar construction can be carried out on the right hand side of the interval.  Replacing $u^{(n)}_{t}$ by  $v_{0} s^{(n)}_{t}$ on the left hand side of the interval and by  $v_{0} s^{(n)}_{t}$ on the right hand side, it is straightforward to check that these unitaries live in $Z_{p,q}$ and implement an approximate unitary equivalence between $\beta$ and $\gamma$.   

(ii), although a little more complicated, follows essentially the same ideas. We omit the proof.
\end{nproof}
\en

\bn
\label{Z-p-q-embedding}
The following is essentially contained in \cite[Theorem~2.1 and Proposition~2.2]{Ror:Z-absorbing}.

\begin{nprop}
Let $\mathfrak{p}$ and $\mathfrak{q}$ be supernatural numbers which are relatively prime. Then, there is a unital embedding 
\[
\bar{\sigma}_{\mathfrak{p},\mathfrak{q}}: Z_{\mathfrak{p},\mathfrak{q}} \to \Zh
\]
which is \emph{standard} in the sense that $\tau_{\Zh} \circ \bar{\sigma}_{\mathfrak{p},\mathfrak{q}} \in T(Z_{\mathfrak{p},\mathfrak{q}})$ is induced by the Lebesgue measure on the unit interval (here, $\tau_{\Zh} \in T(\Zh)$ is the unique tracial state; cf.\ also  (2.1) of \cite{Ror:Z-absorbing}).

If $\sigma: Z_{\mathfrak{p},\mathfrak{q}} \to \Zh$ is another standard embedding, then $\sigma \au \bar{\sigma}_{\mathfrak{p},\mathfrak{q}}$.   
\end{nprop}

\begin{nproof}
\cite[Proposition~2.2]{Ror:Z-absorbing} is only formulated for $\mathfrak{p}$ and $\mathfrak{q}$ of the form $p^{\infty}$ and $q^{\infty}$, respectively, but it is obvious how to adapt the argument to more general relatively prime supernatural numbers. One then obtains increasing sequences $(P_{i})_{i \in \N}$ and $(Q_{i})_{i \in \N}$ of natural numbers and an inductive system
\[
Z_{P_{1},Q_{1}} \to Z_{P_{2},Q_{2}} \to \ldots \to Z_{\mathfrak{p},\mathfrak{q}},
\]
where each connecting map (as well as each limit map $\varrho_{i}$) is a unital $\Ch([0,1])$-homomorphism. As in the proof of \cite[Proposition~2.2]{Ror:Z-absorbing}, there are  standard unital embeddings  $\psi_{i}: Z_{P_{i},Q_{i}} \to \Zh$. Now following the argument of \cite[Proposition~2.2]{Ror:Z-absorbing},  one obtains a unital embedding $\bar{\sigma}_{\mathfrak{p},\mathfrak{q}}: Z_{\mathfrak{p},\mathfrak{q}} \to \Zh$ such that, for each $i$, $\psi_{i} \au \bar{\sigma}_{\mathfrak{p},\mathfrak{q}} \circ \varrho_{i}$. But then it is straightforward to check that  $\bar{\sigma}_{\mathfrak{p},\mathfrak{q}}$ is also standard in the sense of the proposition.

Writing $Z_{\mathfrak{p},\mathfrak{q}}$ as an inductive limit of dimension drop algebras in the obvious way, we obtain the last statement from \cite[Theorem~2.1(ii)]{Ror:Z-absorbing}.
\end{nproof}
\en

\bn
\label{p-Z-p-q-tau-preserving}
\begin{nprop}
Let $\mathfrak{p}$ and $\mathfrak{q}$ be supernatural numbers.  Let 
\[
\varphi:A \otimes Z_{\mathfrak{p},\mathfrak{q}}  \to B \otimes Z_{\mathfrak{p},\mathfrak{q}}
\]
be a $\Ch([0,1])$-homomorphism with fibre maps
\[
\varphi_{t}:(A \otimes Z_{\mathfrak{p},\mathfrak{q}})_{t}  \to (B \otimes Z_{\mathfrak{p},\mathfrak{q}})_{t}
\]
for $t \in [0,1]$. Suppose we have  
\begin{equation}
\label{1c}
\varphi_{t}^{T}(\tau \otimes \bar{\tau}_{t}) = \Delta(\tau) \otimes \bar{\tau}_{t}
\end{equation}
for any $\tau \in T(B)$ and $t \in [0,1]$, where $\bar{\tau}_{t}$ denotes the unique tracial state on the UHF algebra $(Z_{\mathfrak{p},\mathfrak{q}})_{t}$ and $\Delta:T(B) \to T(A)$ is a continuous affine map.  Then, we have 
\[
\varphi^{T}(\tau  \otimes \bar{\tau}_{Z_{\mathfrak{p},\mathfrak{q}}}) = \Delta(\tau ) \otimes \bar{\tau}_{Z_{\mathfrak{p},\mathfrak{q}}} \; \forall \, \tau \in T(B),
\]
where $\bar{\tau}_{Z_{\mathfrak{p},\mathfrak{q}}}$ is the canonical tracial state on $Z_{\mathfrak{p},\mathfrak{q}}$ induced by the Lebesgue measure on $[0,1]$.
\end{nprop}

\begin{nproof}
For any $a \otimes f \in A \otimes Z_{\mathfrak{p},\mathfrak{q}}$, we compute
\begin{eqnarray*}
\varphi^{T}(\tau  \otimes \bar{\tau}_{Z_{\mathfrak{p},\mathfrak{q}}})(a \otimes f) & = & (\tau  \otimes \bar{\tau}_{Z_{\mathfrak{p},\mathfrak{q}}})\circ \varphi (a \otimes f) \\
& = & \int_{[0,1]} (\tau \otimes \bar{\tau}_{t}) (\varphi(a \otimes f)_{t}) \, d\lambda \\
& = & \int_{[0,1]} (\tau \otimes \bar{\tau}_{t}) (\varphi_{t}(a \otimes f_{t})) \, d\lambda \\
& = & \int_{[0,1]} \varphi^{T}_{t}(\tau \otimes \bar{\tau}_{t})(a \otimes f_{t})  \, d\lambda \\
& \stackrel{\eqref{1c}}{=} & \int_{[0,1]} (\Delta(\tau) \otimes \bar{\tau}_{t}) (a \otimes f_{t}) \, d \lambda \\
& = & (\Delta(\tau) \otimes \bar{\tau}_{t}) (a \otimes f).
\end{eqnarray*}
\end{nproof}
\en

\bn
\label{p-phi1-phi2-au}
\begin{nprop}
Let $A$ and $B$ be separable $C^{*}$-algebras, with $B$ unital and $\Zh$-stable; let $\mathfrak{p}$ and $\mathfrak{q}$ be supernatural numbers which are relatively prime. Suppose 
\[
\varphi_{1},\varphi_{2}: A \to B 
\]
are $*$-homomorphisms such that 
\[
\varphi_{1} \otimes \be_{Z_{\mathfrak{p},\mathfrak{q}}} \au \varphi_{2} \otimes \be_{Z_{\mathfrak{p},\mathfrak{q}}}
\]
(as maps from $A$ to $B \otimes Z_{\mathfrak{p},\mathfrak{q}}$). Then,
\[
\varphi_{1} \au \varphi_{2}.
\]
\end{nprop}

\begin{nproof}
Let
\[
\bar{\sigma}_{\mathfrak{p},\mathfrak{q}}: Z_{\mathfrak{p},\mathfrak{q}} \to \Zh
\]
be the embedding from Proposition~\ref{Z-p-q-embedding}. Choose an isomorphism
\[
\nu_{B}:B \stackrel{\cong}{\longrightarrow} B \otimes \Zh
\]
such that
\[
\nu_{B} \au \id_{B} \otimes \be_{\Zh}.
\]
But then, we have
\begin{eqnarray*}
\varphi_{1} & = & \nu_{B}^{-1} \nu_{B} \varphi_{1} \\
& \au & \nu_{B}^{-1} (\id_{B} \otimes \be_{\Zh}) \varphi_{1} \\
& = & \nu_{B}^{-1} (\id_{B} \otimes \bar{\sigma}_{\mathfrak{p},\mathfrak{q}}) (\varphi_{1} \otimes \be_{Z_{\mathfrak{p},\mathfrak{q}}}) \\
& \au & \nu_{B}^{-1}(\id_{B} \otimes \bar{\sigma}_{\mathfrak{p},\mathfrak{q}})  (\varphi_{2} \otimes \be_{Z_{\mathfrak{p},\mathfrak{q}}}) \\
& = & \nu_{B}^{-1} (\id_{B} \otimes \be_{\Zh}) \varphi_{2} \\
& \au & \nu_{B}^{-1} \nu_{B} \varphi_{2} \\
& = & \varphi_{2}.
\end{eqnarray*}
\end{nproof}
\en

\bn
\label{d-exponential-generation}
\begin{ndefn}
Let $\mathfrak{p}$ be a supernatural number of infinite type. Let $(P_{k})_{k\in \N}$ be a sequence of natural numbers. We say $(P_{k})_{k\in \N}$ generates  $\mathfrak{p}$ exponentially, if the following conditions hold:
\begin{enumerate}
\item  $\mathfrak{p}=P_{1}\cdot P_{2} \cdot \ldots$
\item $P_{k}^{2} $ divides $P_{k+1}$ for all $k \in \N$.
\end{enumerate}
\end{ndefn}
\en

\bn
\label{rem-exponential-generation}
\begin{nremarks}
(i) It is easy to see that for any supernatural number of infinite type there is a sequence generating it exponentially. For example, if $\mathfrak{p}=p_{1}^{\infty} \cdot p_{2}^{\infty} \cdot \ldots$, one can set $P_{k}:= p_{1}^{2^{k}} \cdot \ldots \cdot p_{k}^{2^{k}}$.

(ii) If $(P_{k})_{k \in \N}$ generates $\mathfrak{p}$ exponentially, then so does any subsequence $(P_{k_{i}})_{i \in \N}$.

(iii) If  $(P_{k})_{k\in \N}$ generates  $\mathfrak{p}$ exponentially, one has unital inclusions $M_{P_{k}} \subset M_{P_{k+1}}$ and an isomorphism $\lim_{\to} M_{P_{k}} \cong M_{\mathfrak{p}}$.
\end{nremarks}
\en

\bn
\label{not-exponential-generation}
\begin{nnotation}
Let $\mathfrak{p}$ and $\mathfrak{q}$ be  supernatural numbers of infinite type which are relatively prime. Suppose $(P_{k})_{k\in \N}$ and $(Q_{k})_{k\in \N}$ exponentially generate $\mathfrak{p}$ and $\mathfrak{q}$. The inclusions  of Remark~\ref{rem-exponential-generation}(iii) induce inclusions 
\[
\gamma_{k,k+1}: Z_{P_{k},Q_{k}} \hookrightarrow Z_{P_{k+1},Q_{k+1}}.
\]
With this notation, we obtain the inductive system
\[
Z_{\mathfrak{p},\mathfrak{q}}\cong \lim_{\to} (Z_{P_{k},Q_{k}},\gamma_{k,k+1})
\]
with  limit maps
\[
\gamma_{k}: Z_{P_{k},Q_{k}} \hookrightarrow Z_{\mathfrak{p},\mathfrak{q}}. 
\]
For each $k \in \N$, we may choose conditional expectations 
\[
\kappa_{k}: Z_{\mathfrak{p},\mathfrak{q}} \to Z_{P_{k},Q_{k}}.
\]
We shall denote the extensions 
\[
\Ch([0,1]) \otimes M_{\mathfrak{p}} \otimes M_{\mathfrak{q}} \to \Ch([0,1]) \otimes M_{P_{k}} \otimes M_{Q_{k}}
\]
by $\kappa_{k}$ as well.
\end{nnotation}
\en

\bn
\label{def-p-q-exponential-system}
\begin{ndefn}
Let $\mathfrak{p}$ and $\mathfrak{q}$ be supernatural numbers of infinite type which are relatively prime. By a $(\mathfrak{p},\mathfrak{q})$-exponential system for $\Zh$ we mean an inductive system of c.p.c.\ maps
\[
Z_{\mathfrak{p},\mathfrak{q}} \stackrel{\varrho_{1}}{\longrightarrow} Z_{\mathfrak{p},\mathfrak{q}} \stackrel{\varrho_{2}}{\longrightarrow} \ldots 
\]
together with sequences $(P_{k})_{k\in \N} \subset \N$ and $(Q_{k})_{k\in \N} \subset \N$ with the following properties:
\begin{enumerate}
\item[(i)] $\mathfrak{p}$ and $\mathfrak{q}$ are generated exponentially by $(P_{k})_{k\in \N}$ and $(Q_{k})_{k\in \N}$, respectively
\item[(ii)] there are $*$-homomorphisms 
\[
\bar{\varrho}_{k}: Z_{P_{2k},Q_{2k}} \to Z_{P_{2k+1},Q_{2k+1}}
\]
such that  
\[
\varrho_{k}=  \gamma_{2k+1} \circ\bar{\varrho}_{k} \circ \kappa_{2k}
\]
(where $\gamma_{k}$ and $\kappa_{k}$ are as in \ref{not-exponential-generation})
\item[(iii)] for any increasing sequence $(k_{i})_{i \in \N}$, the system $(Z_{\mathfrak{p},\mathfrak{q}},\varrho_{k_{i}})$ is a generalized inductive system (in the sense of \cite{BlaKir:limits}) with limit isomorphic to $\Zh$
\item[(iv)] for any sequence and system as in (iii), there is a sequence of unitaries $(v_{i})_{i \in \N} \subset \Zh$ such that  
\[
\ad(v_{i}) \circ( \ldots \circ \varrho_{k_{i+1}} \circ \varrho_{k_{i}}) \stackrel{i \to \infty}{\longrightarrow} \bar{\sigma}_{\mathfrak{p},\mathfrak{q}} \mbox{ pointwise},
\]
where $\bar{\sigma}_{\mathfrak{p},\mathfrak{q}}$ is as in Proposition~\ref{Z-p-q-embedding}
\item[(v)] for any sequence and system as in (iii), and for any finite subset $\Gh \subset Z_{\mathfrak{p},\mathfrak{q}}$ and $\eta>0$, there is $\bar{m} \in \N$ such that, if $m \ge \bar{m}$, there is a unitary $z \in Z_{P_{2k_{m}+1},Q_{2k_{m}+1}} \otimes Z_{P_{2k_{m}+1},Q_{2k_{m}+1}}$ satisfying
\begin{equation}
\label{6}
\ad (z) \circ (\be_{Z_{P_{2k_{m}+1},Q_{2k_{m}+1}}} \otimes \bar{\varrho}_{k_{m}}) =_{(\Gh,\eta)} \bar{\varrho}_{k_{m}} \otimes \be_{Z_{P_{2k_{m}+1},Q_{2k_{m}+1}}}. 
\end{equation}
\end{enumerate}
\end{ndefn}
\en

\bn
\label{lem-p-q-exponential-system}
\begin{nlemma}
Given supernatural numbers $\mathfrak{p}$ and $\mathfrak{q}$ of infinite type which are relatively prime, there is a  $(\mathfrak{p},\mathfrak{q})$-exponential system for $\Zh$ as in Definition~\ref{def-p-q-exponential-system}. 
\end{nlemma}

\begin{nproof}
By \cite[Theorem~5(a)]{JiaSu:Z},  $M_{\mathfrak{p}}$ and $M_{\mathfrak{q}}$ are both $\Zh$-stable. It then follows from \cite[Theorem~4.6]{HirshbergRordamWinter:D-stable} that $Z_{\mathfrak{p},\mathfrak{q}}$ is also $\Zh$-stable -- in particular, this means that $Z_{\mathfrak{p},\mathfrak{q}}$ admits a unital  embedding $\varrho$ of $\Zh$. Composing this with the embedding $\bar{\sigma}_{\mathfrak{p},\mathfrak{q}}$ from Proposition~\ref{Z-p-q-embedding}, we obtain an endomorphism 
\[
\alpha = \varrho \circ \bar{\sigma}_{\mathfrak{p},\mathfrak{q}}: Z_{\mathfrak{p},\mathfrak{q}} \to Z_{\mathfrak{p},\mathfrak{q}}
\]
with the property that $\alpha^{T}(\tau) \in T(Z_{\mathfrak{p},\mathfrak{q}})$ is induced by the Lebesgue measure on $[0,1]$ for any $\tau \in T(Z_{\mathfrak{p},\mathfrak{q}})$.

Now recall from \ref{not-exponential-generation} that $Z_{\mathfrak{p},\mathfrak{q}}$ can be written as an increasing union of ordinary dimension drop intervals of the form $Z_{P_{k},Q_{k}}$, and  that such algebras are weakly semiprojective by \cite[Proposition~7.3]{JiaSu:Z}.  We then see that $\alpha$ may be approximated pointwise by compositions of maps of the form
\[
\varrho_{k}: Z_{\mathfrak{p},\mathfrak{q}} \stackrel{\kappa_{2k}}{\longrightarrow} Z_{P_{2k},Q_{2k}} \stackrel{\bar{\varrho}_{k}}{\longrightarrow} Z_{P_{2k+1},Q_{2k+1}} \stackrel{\gamma_{2k+1}}{\longrightarrow} Z_{\mathfrak{p},\mathfrak{q}} ,
\]
where the $\bar{\varrho}_{k}$ are  $*$-homomorphisms and the $\kappa_{2k}$ and $\gamma_{2k+1}$ are as in \ref{not-exponential-generation}. It is clear that we may even choose  the sequences $(P_{k})_{k \in \N}$ and $(Q_{k})_{k\in \N}$ to exponentially generate $\mathfrak{p}$ and $\mathfrak{q}$, respectively. Since the $\varrho_{k}$ become more and more multiplicative, with a little extra effort (and possibly passing to subsequences) we may even assume that $(Z_{\mathfrak{p},\mathfrak{q}},\varrho_{k})$ becomes a generalized inductive system.  But then, we clearly have (a composition of)  asymptotic intertwinings (in the sense of \cite{BlaKir:limits})
\begin{equation}
\label{Z-alpha-limit}
\xymatrix{
\ldots \ar[r]  & Z_{P_{2k},Q_{2k}} \ar[r]^{\tilde{\varrho}_{k}} & Z_{P_{2(k+1)},Q_{2(k+1)}} \ar[r] & \ldots \ar[r] & \lim_{\to}(Z_{P_{2k},Q_{2k}},\tilde{\varrho}_{k})\\
\ldots \ar[r] & Z_{\mathfrak{p},\mathfrak{q}} \ar[u]^{\kappa_{2k}} \ar[r]^{\varrho_{k}} & Z_{\mathfrak{p},\mathfrak{q}} \ar[u]^{\kappa_{2(k+1)}} \ar[r] & \ldots \ar[r] & \lim_{\to} (Z_{\mathfrak{p},\mathfrak{q}},\varrho_{k} \ar[u]_{\cong}) \\
\ldots \ar[r] & Z_{\mathfrak{p},\mathfrak{q}} \ar[u]^{=} \ar[r]^{\alpha} & Z_{\mathfrak{p},\mathfrak{q}} \ar[u]^{=} \ar[r] & \ldots \ar[r] & \lim_{\to} (Z_{\mathfrak{p},\mathfrak{q}},\alpha) \ar[u]_{\cong} \ar@/_3.5em/[uu]_{\beta},
}
\end{equation}
where $\tilde{\varrho}_{k} = \kappa_{2(k+1)}  \gamma_{2k+1}\bar{\varrho}_{k}$.  Note that the upper and lower rows are honest inductive systems. Moreover, it is trivial from the pertinent property of $\alpha$ that $\lim_{\to}(Z_{\mathfrak{p},\mathfrak{q}},\alpha)$ has a unique faithful tracial state; from this it is easy to conclude that the limit in fact is a simple $C^{*}$-algebra. But then, $\lim_{\to}(Z_{P_{2k},Q_{2k}},\tilde{\varrho}_{k})$ is a monotracial simple limit of dimension drop intervals, hence may be identified with $\Zh$ by \cite[Theorem~2]{JiaSu:Z}. Therefore, Definition~\ref{def-p-q-exponential-system}(iii) holds. 

We have  established an isomorphism $\beta: \lim_{\to} (Z_{\mathfrak{p},\mathfrak{q}},\alpha) \to \Zh$. Since the limit maps of the lower row of \eqref{Z-alpha-limit} are just given by $\ldots \circ \alpha \circ \alpha$ (starting at stage $k$, say), now Proposition~\ref{Z-p-q-embedding}  shows that  
\[
\beta \circ (\ldots \circ \alpha \circ \alpha) \au \bar{\sigma}_{\mathfrak{p},\mathfrak{q}} .
\]
But then it is clear that, again possibly passing to subsequences, we can even make sure Condition~\ref{def-p-q-exponential-system}(iv) holds. Condition~\ref{def-p-q-exponential-system}(v), possibly oncemore passing to subsequences, now also follows from the inductive limit decomposition \eqref{Z-alpha-limit} of $\Zh$ and the fact that $\Zh$ has approximately inner half flip. 
\end{nproof}
\en

\bn
\begin{nremark}
The proof of the lemma above in particular shows that $\Zh$ can be written as the limit of a stationary inductive system $(Z_{\mathfrak{p},\mathfrak{q}},\alpha)$, where $\alpha$ is any endomorphism of $Z_{\mathfrak{p},\mathfrak{q}}$ with the property that the induced map $\alpha^{T}$ sends all tracial states  on $Z_{\mathfrak{p},\mathfrak{q}}$ to a single one. This point of view will be pursued in more detail in \cite{RordamWinter:Z-revisited}.
\end{nremark}
\en

\section{An asymptotic intertwining}

Below we introduce unitarily suspended isomorphisms between tensor products with generalized dimension drop intervals and show how such an isomorphism may be used to obtain an asymptotic unitary intertwining between tensor products with $\Zh$.

\bn
\label{tensor-factor-not}
\begin{nnotation}
For $0 < m \le n \in \N$ and unital $C^{*}$-algebras $A_{1}, \ldots, A_{m},B_{1},\ldots,B_{n}$ let $i_{1},\ldots,i_{m} \in \{1, \ldots,n$ be pairwise distinct numbers and $\varrho_{j}:A_{j} \to B_{i_{j}}$ be $*$-homomorphisms.  These induce a $*$-homomorphism $$\varrho:=  \varrho_{1} \otimes \ldots \otimes \varrho_{m}: A_{1} \otimes \ldots \otimes A_{m} \to B_{i_{1}} \otimes \ldots \otimes B_{i_{m}}$$.  Composition of this map with the canonical unital embedding
\[
\iota: B_{i_{1}} \otimes \ldots \otimes B_{i_{m}} \to B_{1} \otimes \ldots \otimes B_{n}
\]
will be denoted by $\varrho^{[i_{1}, \ldots,i_{m}]}$, i.e.,
\[
\varrho^{[i_{1}, \ldots,i_{m}]} :=  \iota  \circ  \varrho.
\]
With this notation, $\iota $ above may be expressed as $\id^{[i_{1}, \ldots,i_{m}]}$. Similarly, the flip on a $C^{*}$-algebra $\Dh$ can be written as $\id^{[2,1]}$.

For $1 \le k \le m$, we will sometimes have unital $*$-homomorphisms of the form 
\[
\sigma: \bigotimes_{j \in \{1, \ldots \check{k} \ldots,m} A_{j} \to B_{i_{j}}.
\]
The maps $\sigma$ and $\varrho_{i_{k}}$ then give rise to a unital $*$-homomorphism 
\[
\sigma^{[i_{1}, \ldots \check{i_{k}} \ldots,i_{m}]} \otimes \varrho^{[i_{k}]}: A_{1} \otimes \ldots \otimes A_{m} \to B_{1} \otimes \ldots \otimes B_{n} 
\]
in the obvious way.
\end{nnotation}
\en

\bn
\label{d-unitarily-suspended}
\begin{ndefn}
Let $A$ and $B$ be unital $C^{*}$-algebras and let $\mathfrak{p}$ and $\mathfrak{q}$ be supernatural numbers. We say a $\Ch([0,1])$-homomorphism 
\[
\varphi: A \otimes Z_{\mathfrak{p},\mathfrak{q}} \to B \otimes Z_{\mathfrak{p},\mathfrak{q}}
\]
is unitarily suspended, if there are $0 \le t_{0} < t_{1} \le 1$, a continuous path $(u_{t})_{t \in [t_{0},t_{1})}$ of unitaries in $B \otimes M_{\mathfrak{p}} \otimes M_{\mathfrak{q}}$ and   $*$-homomorphisms
\[
\sigma_{\mathfrak{p}}: A \otimes M_{\mathfrak{p}} \to B \otimes M_{\mathfrak{p}}
\]
and 
\[
\varrho_{\mathfrak{q}}: A \otimes M_{\mathfrak{q}} \to B \otimes M_{\mathfrak{q}}
\]
such that $u_{t_{0}}=\be_{B \otimes M_{\mathfrak{p}} \otimes M_{\mathfrak{q}}}$ and
\begin{equation}
\label{77}
\varphi^{(t)}= \left\{
\begin{array}{lll}
\sigma_{\mathfrak{p}} & \mbox{for} & t = 0\\
\sigma_{\mathfrak{p}} \otimes \id_{M_{\mathfrak{q}}} & \mbox{for} & t \in (0,t_{0}] \\
\ad (u_{t}) \circ (\sigma_{\mathfrak{p}} \otimes \id_{M_{\mathfrak{q}}}) & \mbox{for} & t \in (t_{0},t_{1}) \\
\varrho_{\mathfrak{q}}^{[1,3]} \otimes \id_{M_{\mathfrak{p}}}^{[2]} & \mbox{for} & t \in [t_{1},1) \\
\varrho_{\mathfrak{q}} & \mbox{for} & t = 1.
\end{array}
\right.
\end{equation}
\end{ndefn}
\en

\bn
\label{rho-m-rho-n-intertwining}
The following Lemma will be used in Proposition~\ref{barbarphi-intertwining} to construct a certain approximate unitary intertwining. The two results together mark the crucial technical step in the proof of Theorem~\ref{main-result}.

\begin{nlemma}
Let $\mathfrak{p}$ and $\mathfrak{q}$ be supernatural numbers of infinite type which are relatively prime. 
Suppose $A$ and $B$ are separable, unital and $\Zh$-stable $C^{*}$-algebras; suppose furthermore that  
\[
\varphi: A \otimes Z_{\mathfrak{p},\mathfrak{q}} \to B  \otimes Z_{\mathfrak{p},\mathfrak{q}} 
\]
is a unitarily suspended $\Ch([0,1])$-isomorphism.

Then, there is a $(\mathfrak{p},\mathfrak{q})$-exponential system for $\Zh$ as in Definition~\ref{def-p-q-exponential-system} with the following property:  

Given  a finite subset $\Fh \subset A  \otimes Z_{\mathfrak{p},\mathfrak{q}}$ and  $\varepsilon>0$, there is $\bar{m} \in \N$ such that the following holds: for any $m \ge \bar{m}$ there is $\bar{n} \in \N$ such that, if $n \ge \bar{n}$, there is a unitary $u \in B  \otimes Z_{\mathfrak{p},\mathfrak{q}}$ satisfying
\[
\varphi \circ (\id_{A} \otimes \varrho_{m}) =_{(\Fh,\epsilon)} \ad (u) \circ (\id_{B} \otimes \varrho_{n}) \circ \varphi \circ (\id_{A} \otimes \varrho_{m}).
\]
\end{nlemma}

\begin{nproof}
Let
\begin{equation}
\label{8c}
\xymatrix{
\ldots \ar[r] & Z_{\mathfrak{p},\mathfrak{q}} \ar@/^2em/[rrr]^{\varrho_{k}} \ar[r]^{\kappa_{2k}} & Z_{P_{2k},Q_{2k}} \ar[r]^{\bar{\varrho}_{k}} & Z_{P_{2k+1},Q_{2k+1}} \ar[r]^{\gamma_{2k+1}} & Z_{\mathfrak{p},\mathfrak{q}} \ar[r] & \ldots 
}
\end{equation}
be a $(\mathfrak{p},\mathfrak{q})$-exponential system for $\Zh$; such a system exists by Lemma~\ref{lem-p-q-exponential-system}.

For convenience, we  set 
\begin{equation}
\label{19}
\eta:= \frac{\varepsilon}{12}.
\end{equation}

Since $A$ and $B$ are $\Zh$-stable, we may as well  replace them by $A \otimes \Zh$ and $B \otimes \Zh$, respectively. Using the notation of Definition~\ref{d-unitarily-suspended} for the unitarily suspended map $\varphi: A \otimes \Zh \otimes Z_{\mathfrak{p},\mathfrak{q}} \to A \otimes \Zh \otimes Z_{\mathfrak{p},\mathfrak{q}}$, we set
\begin{equation}
\label{3c}
\sigma:= \sigma_{\mathfrak{p}} \otimes \id_{M_{\mathfrak{q}}} : A \otimes \Zh \otimes M_{\mathfrak{p}} \otimes M_{\mathfrak{q}} \to B \otimes \Zh \otimes M_{\mathfrak{p}} \otimes M_{\mathfrak{q}}.
\end{equation}

We may assume $\Fh$ to be of the form 
\[
\Fh= \{ a \otimes c \mid a \in \Fh',\, c \in \Fh''\}
\]
for  finite subsets of normalized elements $$\be_{A \otimes \Zh} \in \Fh' \subset (A \otimes \Zh)_{+}$$ and $$\be_{Z_{P_{2k_{m'}},Q_{2k_{m'}}}} \in \Fh'' \subset \gamma_{2k_{m'}}  (Z_{P_{2k_{m'}},Q_{2k_{m'}}})_{+} \subset (Z_{\mathfrak{p},\mathfrak{q}})_{+},$$ where $k_{m'} \in \N$ is some natural number. 

By \cite[Remark~2.7]{TomsWinter:ssa}, we may furthermore assume that there is an isomorphism 
\[
\nu_{A}: A \otimes \Zh \otimes \Zh \stackrel{\cong}{\longrightarrow} A \otimes \Zh
\]
such that
\[
\nu_{A} \circ (\id_{A \otimes \Zh} \otimes \be_{\Zh}) =_{(\Fh',\eta)} \id_{A \otimes \Zh}
\]
and
\[
\nu_{A} \circ (\id_{A \otimes \Zh} \otimes \be_{\Zh}) \au \id_{A \otimes \Zh}.
\]
Now by Definition~\ref{def-p-q-exponential-system}(v), there is $k_{\bar{m}}\ge k_{m'} \in \N$ such that for any given $k_{m} \ge k_{\bar{m}}$ there is a unitary 
\begin{equation}
\label{6a}
z \in Z_{P_{2k_{m}+1},Q_{2k_{m}+1}} \otimes Z_{P_{2k_{m}+1},Q_{2k_{m}+1}}
\end{equation}
satisfying
\begin{equation}
\label{6b}
\ad (z) \circ (\be_{Z_{P_{2k_{m}+1},Q_{2k_{m}+1}}} \otimes \bar{\varrho}_{k_{m}} \circ \kappa_{2k_{m}}) =_{(\Fh'',\eta)} \bar{\varrho}_{k_{m}} \circ \kappa_{2k_{m}} \otimes \be_{Z_{P_{2k_{m}+1},Q_{2k_{m}+1}}}. 
\end{equation}

For each $k$, we set
\[
C_{k}:= Z_{P_{k},Q_{k}}.
\]
Let 
\[
\theta : C_{2k_{m}+1} \to \Zh
\]
be some unital embedding.

Since
\[
\gamma_{2k_{n}} \circ \kappa_{2k_{n}} \stackrel{n \to \infty}{\longrightarrow} \id_{Z_{\mathfrak{p},\mathfrak{q}}},
\]
and since 
\[
M_{\mathfrak{p}} \otimes M_{\mathfrak{q}} = \lim_{\to} M_{P_{2k_{n}}} \otimes M_{Q_{2k_{n}}},
\]
there is $k_{\bar{n}} \in \N$ such that, if $k_{n} \ge k_{\bar{n}}$ is given, we have 
\begin{equation}
\label{8}
(\id_{B \otimes \Zh} \otimes \gamma_{2k_{n}} \circ \kappa_{2k_{n}}) \circ \varphi \circ (\id_{A \otimes \Zh} \otimes \varrho_{k_{m}}) =_{(\Fh,\eta)} \varphi \circ (\id_{A \otimes \Zh} \otimes \varrho_{k_{m}})
\end{equation}
and, using \ref{d-unitarily-suspended},  a unitary 
\begin{equation}
\label{2c}
u \in B \otimes \Zh \otimes \Ch([0,1], M_{P_{2k_{n}}} {\otimes} M_{Q_{2k_{n}}}) (\cong B \otimes \Zh \otimes \Ch([0,1]) \otimes M_{P_{2k_{n}}} \otimes M_{Q_{2k_{n}}})
\end{equation}
such that 
\[
u(0)=1
\]
and
\begin{equation}
\label{8a}
\ad(u) \circ \sigma^{[1,2,4,5]} \circ (\id_{A\otimes \Zh} \otimes \varrho_{k_{m}}) =_{(\Fh,\eta)} \varphi \circ (\id_{A \otimes \Zh} \otimes \varrho_{k_{m}}),
\end{equation}
where we have identified $B \otimes \Zh \otimes \Ch([0,1], M_{P_{2k_{n}}} \otimes M_{Q_{2k_{n}}}$ with $B \otimes \Zh \otimes \Ch([0,1]) \otimes M_{P_{2k_{n}}} \otimes M_{Q_{2k_{n}}} \subset B \otimes \Zh \otimes \Ch([0,1]) \otimes M_{\mathfrak{p}} \otimes M_{\mathfrak{q}}$.

Moreover, we may assume that there is a unitary 
\begin{equation}
\label{11}
u' \in B \otimes \Zh \otimes \Zh \otimes Z_{\mathfrak{p},\mathfrak{q}}
\end{equation}
with 
\begin{equation}
\label{12}
\|u' - \tilde{u}\| < \eta,
\end{equation}
where
\begin{equation}
\label{13}
\tilde{u} := (((\id_{B \otimes \Zh} \otimes \varrho_{k_{n}})^{[1,2,4]} \circ \varphi \circ (\be_{A \otimes \Zh} \otimes \gamma_{2k_{m}+1})) \otimes \theta^{[3]})(z)
\end{equation}
is close to a unitary if only $n$ is large enough because $\varrho_{k_{n}}$ is almost multiplicative (here, $z$ comes from \eqref{6a}). 

Finally, for large enough $n$,  we may even assume
\begin{eqnarray}
\label{9}
\lefteqn{
\ad(\tilde{u}) \circ (((\id_{B \otimes \Zh} \otimes \varrho_{k_{n}})^{[1,2,4]} \circ \varphi \circ (\id_{A \otimes \Zh} \otimes \be_{Z_{\mathfrak{p},\mathfrak{q}}})) \otimes  (\theta \circ \bar{\varrho}_{k_{m}} \circ \kappa_{2k_{m}})^{[3]}) (a \otimes c) }\nonumber \\
& =_{\eta} &  
 (((\id_{B \otimes \Zh} \otimes \varrho_{k_{n}})^{[1,2,4]} \circ \varphi) \otimes \id_{\Zh}^{[3]}) \nonumber \\
 & & \circ \ad((\be_{A \otimes \Zh} \otimes \gamma_{2k_{m}+1} \otimes \theta)(z))(a \otimes \gamma_{2k_{m}+1}(\be_{C_{2k_{m}+1}}) \otimes \theta \circ \bar{\varrho}_{k_{m}} \circ \kappa_{2k_{m}}^{[3]}(c)) 
\end{eqnarray}
for $a \otimes c \in \Fh$ and $n \ge \bar{n}$. 

For the remainder of the proof, fix such an $n \in \N$.

For $l\in \N$, choose $\Ch([0,1])$-maps $\bar{\gamma}_{2k_{l}}$, $\bar{\gamma}_{2k_{l}+1}$ such that the diagram
\begin{equation}
\label{7c}
\xymatrix{
\ldots \ar[r] & C_{2k_{l}} \ar[r]^{\bar{\gamma}_{2k_{l}}} \ar[drrr]_{\gamma_{2k_{l}}} & C_{2k_{l}+1} \ar[r]^{\bar{\gamma}_{2k_{l}+1}} \ar[drr]^{\gamma_{2k_{l}+1}} & C_{2k_{l+1}} \ar[r] \ar[dr]^{\gamma_{2k_{l+1}}} & \ldots \\
&&&& Z_{\mathfrak{p},\mathfrak{q}}
}
\end{equation}
commutes. Using property~\ref{d-exponential-generation}(ii), it is easy to find a $\Ch([0,1])$-map 
\[
\bar{\beta}: C_{2k_{n}} \to C_{2k_{n+1}}
\]
such that
\begin{equation}
\label{5a}
[\bar{\beta}(C_{2k_{n}}),\bar{\gamma}_{2k_{n}+1} \circ \bar{\varrho}_{k_{n}}(C_{2k_{n}})] = 0.
\end{equation}
By Proposition~\ref{p-dimension-drop-inclusions-au}a), we have 
\begin{equation}
\label{4}
\bar{\beta} \au \bar{\gamma}_{2k_{n}+1} \circ \bar{\gamma}_{2k_{n}}.
\end{equation}
Set
\begin{equation}
\label{6c}
\beta:= \bar{\beta} \circ \kappa_{2k_{n}} : Z_{\mathfrak{p},\mathfrak{q}} \to C_{2k_{n+1}}.
\end{equation}
By \eqref{5a}, we have a $*$-homomorphism
\begin{equation}
\label{5c}
\bar{\beta} \otimes (\bar{\gamma}_{2k_{n}+1} \circ \bar{\varrho}_{k_{n}}): C_{2k_{n}} \otimes C_{2k_{n}} \to C_{2k_{n+1}}
\end{equation}
given by
\[
\bar{\beta} \otimes (\bar{\gamma}_{2k_{n}+1} \circ \bar{\varrho}_{k_{n}})(e \otimes f)= \bar{\beta} (e)  (\bar{\gamma}_{2k_{n}+1} \circ \bar{\varrho}_{k_{n}})(f).
\]

Define 
\begin{eqnarray*}
D & := & \{f \in \Ch([0,1]^{2},M_{P_{2k_{n}}} \otimes M_{Q_{2k_{n}}}) \mid \\
&& f(0,0) \in M_{P_{2k_{n}}} \otimes \be_{Q_{2k_{n}}}, \, f(1,1) \in \be_{P_{2k_{n}}} \otimes M_{Q_{2k_{n}}}\}
\end{eqnarray*}
and embeddings
\[
\zeta_{x},\zeta_{y}: C_{2k_{n}} \to D 
\]
as in Proposition~\ref{p-dimension-drop-inclusions-au}(ii) by
\begin{equation}
\label{8b}
\zeta_{x}(f)(s,t) := f(s) \mbox{ and } \zeta_{y}(f)(s,t):= f(t). 
\end{equation}
By Proposition~\ref{p-dimension-drop-inclusions-au}(ii), there is a $\Ch([0,1]^{2})$-embedding
\[
\alpha: D \to C_{2k_{n}} \otimes C_{2k_{n}}
\]
such that  
\begin{equation}
\label{1}
\alpha \circ \zeta_{x} \au \id_{C_{2k_{n}}} \otimes \be_{C_{2k_{n}}}
 \mbox{ and } 
\alpha \circ \zeta_{y} \au \be_{C_{2k_{n}}} \otimes \id_{C_{2k_{n}}}.
\end{equation}

Next, we define a unitary
\[
w \in B \otimes \Zh \otimes D
\]
by
\[
w(s,t):= u(t) u^{*}(s)
\]
(cf.\ \eqref{2c} and \eqref{8a} for the definition of $u$). 

Now for any $(s,t) \in [0,1]^{2} \setminus \{(0,0),(1,1)\}$  and $a \in \Fh'$ we have
\begin{eqnarray*}
\lefteqn{
((\id_{B \otimes \Zh} \otimes \zeta_{y} \circ \kappa_{2k_{n}}) \circ \varphi \circ (\id_{A \otimes \Zh} \otimes \be_{Z_{\mathfrak{p},\mathfrak{q}}})(a))(s,t)
}\\
& = & ((\id_{B \otimes \Zh} \otimes \kappa_{2k_{n}}) \circ \varphi \circ (\id_{A \otimes \Zh} \otimes \be_{Z_{\mathfrak{p},\mathfrak{q}}})(a))(s,t)  \\
& \stackrel{\eqref{8a}}{=}_{\eta} & (\ad(u) \circ (\id_{B \otimes \Zh} \otimes \kappa_{2k_{n}}) \circ \sigma^{[1,2,4,5]} \circ (\id_{A \otimes \Zh} \otimes \be_{Z_{\mathfrak{p},\mathfrak{q}}})(a) )(t)\\
& = & (\ad(u) \circ (\id_{B \otimes \Zh} \otimes \kappa_{2k_{n}}) \circ  (\sigma^{[1,2,4,5]} \otimes \id^{[3]}_{\Ch([0,1])}) \\
& & \circ (\id_{A \otimes \Zh}^{[1,2]} \otimes \varrho_{2k_{m}}^{[3,4,5]})(a \otimes \be_{Z_{\mathfrak{p},\mathfrak{q}}}))(t) \\
& = & (\ad(u) \circ (\id_{B \otimes \Zh} \otimes \kappa_{2k_{n}}) \circ  \sigma^{[1,2,4,5]}(a \otimes \be_{M_{P_{2k_{m}}}} \otimes \be_{M_{Q_{2k_{m}}}}))(t) \\
& = & \ad(u(t)) ((\id_{B \otimes \Zh} \otimes \kappa_{2k_{n}}) \circ \sigma^{[1,2,4,5]} (a \otimes \be_{M_{P_{2k_{m}}}} \otimes \be_{M_{Q_{2k_{m}}}})(s))\\
& = & \ad(u(t) u^{*}(s)u(s)) \sigma^{[1,2,4,5]}(a \otimes \be_{M_{P_{2k_{m}}}}\otimes \be_{M_{Q_{2k_{m}}}})(s) \\
& \vdots & \\
& =_{\eta} & \ad(w(s,t))(((\id_{B \otimes \Zh} \otimes \zeta_{x} \circ \kappa_{2k_{n}}) \circ \varphi \circ (\id_{A \otimes \Zh} \otimes \be_{Z_{\mathfrak{p},\mathfrak{q}}})(a))(s,t)) \\
& = & (\ad(w) \circ (\id_{B \otimes \Zh} \otimes \zeta_{x} \circ \kappa_{2k_{n}}) \circ \varphi \circ (\id_{A \otimes \Zh} \otimes \be_{Z_{\mathfrak{p},\mathfrak{q}}})(a))(s,t), \\
\end{eqnarray*}
where for the second part we argue in essentially the same way as for the first part. For $(s,t)=(0,0)$ or $(s,t)=(1,1)$ the above estimate is trivial, since $w(0,0)= \be_{D}(0,0)$ and $w(1,1)= \be_{D}(1,1)$, and by \eqref{8b}. We obtain
\begin{eqnarray}
\label{2}
\lefteqn{(\id_{B \otimes \Zh} \otimes \zeta_{y} \circ \kappa_{2k_{n}}) \circ \varphi \circ (\id_{A \otimes \Zh} \otimes \be_{Z_{\mathfrak{p},\mathfrak{q}}})} \nonumber \\
& =_{(\Fh',2\eta)} & \ad(w) \circ (\id_{B \otimes \Zh} \otimes \zeta_{x} \circ \kappa_{2k_{n}}) \circ \varphi \circ (\id_{A \otimes \Zh} \otimes \be_{Z_{\mathfrak{p},\mathfrak{q}}}).
\end{eqnarray}
Applying $\id_{B \otimes \Zh} \otimes \alpha$ on both sides of \eqref{2}, and using \eqref{1}, we see that there is a unitary
\[
w' \in B \otimes \Zh \otimes C_{2k_{n}} \otimes C_{2k_{n}}
\]
satisfying
\begin{eqnarray}
\label{3}
\lefteqn{(\id_{B \otimes \Zh} \otimes \be_{C_{2k_{n}}} \otimes \kappa_{2k_{n}}) \circ \varphi \circ (\id_{A \otimes \Zh} \otimes \be_{Z_{\mathfrak{p},\mathfrak{q}}})} \nonumber \\
& =_{(\Fh',3\eta)} & \ad(w') \circ (\id_{B \otimes \Zh} \otimes \kappa_{2k_{n}} \otimes \be_{C_{2k_{n}}}) \circ \varphi \circ (\id_{A \otimes \Zh} \otimes \be_{Z_{\mathfrak{p},\mathfrak{q}}}).
\end{eqnarray}
But then,
\[
w'':= (\id_{B \otimes \Zh} \otimes \bar{\beta} \otimes (\bar{\gamma}_{2k_{n}+1} \circ \bar{\varrho}_{k_{n}}))(w') \in B \otimes \Zh \otimes C_{2k_{n+1}}
\]
(cf.\ \eqref{5c}) satisfies
\begin{eqnarray}
\label{10}
\lefteqn{(\id_{B \otimes \Zh} \otimes (\bar{\gamma}_{2k_{n}+1} \circ \bar{\varrho}_{k_{n}} \circ \kappa_{2k_{n}})) \circ \varphi \circ (\id_{A \otimes \Zh} \otimes \be_{Z_{\mathfrak{p},\mathfrak{q}}})} \nonumber \\
& \stackrel{\eqref{5c},\eqref{6c}}{=}_{(\Fh',3\eta)} & \ad(w'') \circ (\id_{B \otimes \Zh} \otimes \beta) \circ \varphi \circ (\id_{A \otimes \Zh} \otimes \be_{Z_{\mathfrak{p},\mathfrak{q}}}).
\end{eqnarray}
Now from \eqref{8c}, \eqref{7c}, \eqref{4} and \eqref{8}, and applying $\id_{B \otimes \Zh} \otimes \gamma_{2k_{n+1}}$ to \eqref{10}, we obtain a unitary
\[
w''' \in B \otimes \Zh \otimes Z_{\mathfrak{p},\mathfrak{q}}
\]
satisfying
\begin{eqnarray}
\label{18}
\lefteqn{(\id_{B \otimes \Zh} \otimes \varrho_{k_{n}}) \circ \varphi \circ (\id_{A \otimes \Zh} \otimes \be_{Z_{\mathfrak{p},\mathfrak{q}}})} \nonumber\\
& \stackrel{\eqref{4},\eqref{6c},\eqref{10}}{=}_{(\Fh',4 \eta)} & \ad(w''') \circ (\id_{B \otimes \Zh} \otimes (\gamma_{2k_{n}} \circ \kappa_{2k_{n}})) \circ \varphi \circ (\id_{A \otimes \Zh} \otimes \be_{Z_{\mathfrak{p},\mathfrak{q}}}) \nonumber \\
& \stackrel{\eqref{8}}{=}_{(\Fh',\eta)} & \ad(w''') \circ \varphi \circ (\id_{A \otimes \Zh} \otimes \be_{Z_{\mathfrak{p},\mathfrak{q}}}).
\end{eqnarray}

Next, employ \cite[Remark~2.7]{TomsWinter:ssa} to choose
\[
\nu_{B}: B \otimes \Zh \otimes \Zh \stackrel{\cong}{\longrightarrow} B \otimes \Zh
\]
such that
\begin{eqnarray}
\label{8cc}
\lefteqn{(\nu_{B} \otimes \id_{Z_{\mathfrak{p},\mathfrak{q}}}) \circ ((\id_{B \otimes \Zh} \otimes \varrho_{k_{n}})\circ \varphi)^{[1,2,4]} \circ (\id_{A \otimes \Zh} \otimes \varrho_{k_{m}})} \nonumber \\
& =_{(\Fh,\eta)} & (\id_{B \otimes \Zh} \otimes \varrho_{k_{n}}) \circ \varphi \circ (\id_{A \otimes \Zh} \otimes \varrho_{k_{m}})
\end{eqnarray}
and such that
\begin{equation}
\label{14}
 \varphi \circ (\id_{A \otimes \Zh} \otimes \varrho_{k_{m}})
 =_{\Fh,\eta}   (\nu_{B} \otimes \id_{Z_{\mathfrak{p},\mathfrak{q}}})\circ \varphi^{[1,2,4]} \circ (\id_{A \otimes \Zh} \otimes \varrho_{k_{m}}). 
\end{equation}
We set
\begin{equation}
\label{17}
u'' := \id_{B \otimes \Zh \otimes Z_{\mathfrak{p},\mathfrak{q}}}^{[1,2,4]}(w''')
\end{equation}
and
\begin{equation}
\label{15}
u''':= (\varphi^{[1,2,4]} \otimes \id_{\Zh}^{[3]}) \circ (\theta \otimes \gamma_{2k_{m}+1})^{[3,4]}(z)
\end{equation}
and define
\begin{equation}
\label{35}
u:= \nu_{B} \otimes \id_{Z_{\mathfrak{p},\mathfrak{q}}}(u' u'' u''')^{*}
\end{equation}
(see \eqref{11} and \eqref{12} for the definition of $u'$).

For $a \in \Fh'$ and $c \in \Fh''$ we are finally ready to compute for $a \in \Fh'$, $c \in \Fh''$
\begin{eqnarray*}
\lefteqn{\ad(u^{*}) \circ \varphi \circ (\id_{A \otimes \Zh} \otimes \varrho_{k_{m}})(a \otimes c)}\\
& \stackrel{\eqref{14}}{=_{\eta}} & \ad(u^{*}) \circ (\nu_{B} \otimes \id_{Z_{\mathfrak{p},\mathfrak{q}}})\circ \varphi^{[1,2,4]} \circ (\id_{A \otimes \Zh} \otimes \varrho_{k_{m}})(a \otimes c) \\
& \stackrel{\eqref{15},\eqref{6b},\eqref{35}}{=_{\eta}} & (\nu_{B} \otimes \id_{Z_{\mathfrak{p},\mathfrak{q}}}) \circ \ad(u' u'') \circ (\varphi^{[1,2,4]}(a \otimes \be_{Z_{\mathfrak{p},\mathfrak{q}}}) \otimes (\theta \circ \bar{\varrho}_{k_{m}} \circ \kappa_{2k_{m}})^{[3]}(c)) \\
& \stackrel{\eqref{17},\eqref{18}}{=_{5\eta}} & (\nu_{B} \otimes \id_{Z_{\mathfrak{p},\mathfrak{q}}}) \circ \ad(u') \\
& & \circ (((\id_{B \otimes \Zh} \otimes \varrho_{k_{n}}) \circ \varphi)^{[1,2,4]}(a \otimes \be_{Z_{\mathfrak{p},\mathfrak{q}}}) \otimes (\theta \circ \bar{\varrho}_{k_{m}} \circ \kappa_{2k_{m}})^{[3]}(c)) \\
& \stackrel{\eqref{12}}{=_{2 \eta}} & (\nu_{B} \otimes \id_{Z_{\mathfrak{p},\mathfrak{q}}})  \circ \ad(\tilde{u}) \\
& & \circ (((\id_{B \otimes \Zh} \otimes \varrho_{k_{n}}) \circ \varphi)^{[1,2,4]}(a \otimes \be_{Z_{\mathfrak{p},\mathfrak{q}}}) \otimes (\theta \circ \bar{\varrho}_{k_{m}} \circ \kappa_{2k_{m}})^{[3]}(c)) \\ 
& \stackrel{\eqref{9}}{=_{\eta}} &  (\nu_{B} \otimes \id_{Z_{\mathfrak{p},\mathfrak{q}}})   \circ  (((\id_{B \otimes \Zh} \otimes \varrho_{k_{n}}) \circ \varphi)^{[1,2,4]} \otimes \id_{\Zh}^{[3]}) \\
& & \circ \ad((\be_{A \otimes \Zh} \otimes \gamma_{2k_{m}+1} \otimes \theta)(z))(a \otimes \gamma_{2k_{m}+1}(\be_{C_{2k_{m}+1}}) \otimes \theta  \bar{\varrho}_{k_{m}}  \kappa_{2k_{m}}(c)) \\
& \stackrel{\eqref{6b}}{=_{\eta}} & (\nu_{B} \otimes \id_{Z_{\mathfrak{p},\mathfrak{q}}})   \\
& &\circ  (((\id_{B \otimes \Zh} \otimes \varrho_{k_{n}}) \circ \varphi)^{[1,2,4]} \otimes \id_{\Zh}^{[3]}) (a \otimes \gamma_{2k_{m}+1}  \bar{\varrho}_{k_{m}}  \kappa_{2k_{m}}(c) \otimes \be_{\Zh}) \\
& = & (\nu_{B} \otimes \id_{Z_{\mathfrak{p},\mathfrak{q}}})   \circ  (((\id_{B \otimes \Zh} \otimes \varrho_{k_{n}}) \circ \varphi)^{[1,2,4]}) (a \otimes \varrho_{k_{m}}(c)) \\
& \stackrel{\eqref{8cc}}{=_{\eta}} & (\id_{B \otimes \Zh} \otimes \varrho_{k_{n}}) \circ \varphi \circ (\id_{A \otimes \Zh} \otimes \varrho_{k_{m}})(a \otimes c).
\end{eqnarray*}
Therefore (and by \eqref{19}) we have
\begin{equation}
\label{9c}
\varphi \circ (\id_{A \otimes \Zh} \otimes \varrho_{k_{m}}) =_{(\Fh,\varepsilon)} \ad(u) \circ (\id_{B \otimes \Zh} \otimes \varrho_{k_{n}}) \circ \varphi \circ (\id_{A \otimes \Zh} \otimes \varrho_{k_{m}}).
\end{equation}
We may now pass to a subsystem of \eqref{8c}; by relabeling the indices, we obtain the assertion of the lemma from \eqref{9c}.
\end{nproof}
\en

\bn
\label{barbarphi-intertwining}
\begin{nprop}
Assume the same hypotheses as in \ref{rho-m-rho-n-intertwining} and fix a  $(\mathfrak{p},\mathfrak{q})$-exponential system for $\Zh$ such that the assertion of  \ref{rho-m-rho-n-intertwining} holds. 

Then, there are increasing sequences $(m_{i})_{i \in \N}$, $(n_{i})_{i \in \N} \subset \N$  which yield a  diagram
\begin{equation}
\label{barbarphi-diagram}
\xymatrix @C=-0.7em @W=3.5em @R7ex{
\cdots \ar[rr]& & A {\otimes} Z_{\mathfrak{p},\mathfrak{q}} \ar[dr]_(.4){\alpha_{i}} \ar[rr]^{\id_{A} \otimes \varrho_{m_{i}}} && A {\otimes} Z_{\mathfrak{p},\mathfrak{q}} \ar[dr] \ar[rr] &&\cdots &   \ar[rr]& & A {\otimes} \Zh \ar@<-0.5em>[d]_{\bar{\bar{\varphi}}} \\
&\cdots  \ar[rr] && B {\otimes} Z_{\mathfrak{p},\mathfrak{q}} \ar[ur]^(.4){\beta_{i}} \ar[rr]^{\id_{B} \otimes \varrho_{n_{i}}} && B  {\otimes} Z_{\mathfrak{p},\mathfrak{q}} \ar[ur] \ar[rr] &&    \cdots \ar[rr] & &B  {\otimes} \Zh \ar@<-0.5em>[u]_{\bar{\bar{\varphi}}^{-1}}
}
\end{equation}
where
\begin{equation}
\label{31}
\alpha_{i}:= \varphi \circ (\id_{A} {\otimes} \varrho_{m_{i}})
\end{equation}
and 
\begin{equation}
\label{32}
\beta_{i}:= \varphi^{-1} \circ (\id_{B} {\otimes} \varrho_{n_{i}}),
\end{equation}
which satisfies the following properties:
\begin{enumerate}
\item[(i)] the diagram \eqref{barbarphi-diagram} is an approximate intertwining in the sense that there are sequences 
\[
(w_{i}^{A})_{i\in \N} \subset A \otimes \Zh
\]
and
\[
(w_{i}^{B})_{i \in \N} \subset B \otimes \Zh
\]
of unitaries such that  the maps $\bar{\bar{\varphi}}$ and $\bar{\bar{\varphi}}^{-1}$,  given by
\begin{eqnarray*}
\lefteqn{\bar{\bar{\varphi}} \circ (\id_{A} \otimes ( \ldots \circ \varrho_{m_{j+2}} \circ \varrho_{m_{j+1}} \circ \varrho_{m_{j}}))}\\
& = & \lim_{i \to \infty} \ad(w_{i}^{B}) \circ ( \id_{B} \otimes (\ldots \circ \varrho_{n_{i+1}} \circ \varrho_{n_{i}})) \circ \varphi \circ (\id_{A} \otimes (\varrho_{m_{i}} \circ \ldots \circ \varrho_{m_{j}})), 
\end{eqnarray*}
and
\begin{eqnarray*}
\lefteqn{\bar{\bar{\varphi}}^{-1} \circ (\id_{B} \otimes ( \ldots \circ \varrho_{n_{j+2}} \circ \varrho_{n_{j+1}} \circ \varrho_{n_{j}}))}\\
& = & \lim_{i \to \infty} \ad(w_{i}^{A}) \circ ( \id_{A} \otimes (\ldots \circ \varrho_{m_{i+1}} \circ \varrho_{m_{i}})) \circ \varphi^{-1} \circ (\id_{B} \otimes (\varrho_{n_{i}} \circ \ldots \circ \varrho_{n_{j}})), 
\end{eqnarray*}
are well-defined and mutually inverse $*$-isomorphisms
\item[(ii)] the map $\bar{\bar{\varphi}}$ satisfies
\[
\bar{\bar{\varphi}} \circ (\id_{A} \otimes \be_{\Zh}) \au (\id_{B} \otimes \bar{\sigma}_{\mathfrak{p},\mathfrak{q}}) \circ \varphi \circ (\id_{A} \otimes \be_{Z_{\mathfrak{p},\mathfrak{q}}})
\]
(and the respective statement holds for $\bar{\bar{\varphi}}^{-1}$), where $\bar{\sigma}_{\mathfrak{p},\mathfrak{q}}: Z_{\mathfrak{p},\mathfrak{q}} \to \Zh$ is the standard embedding from Proposition~\ref{Z-p-q-embedding}.
\end{enumerate}
\end{nprop}

\begin{nproof}
Let $(a_{l})_{l \in \N} \subset A \otimes Z_{\mathfrak{p},\mathfrak{q}}$ and $(b_{l})_{l\in \N} \subset B \otimes Z_{\mathfrak{p},\mathfrak{q}}$ be dense subsets. Set $$\Gh_{0}^{A}:= \{a_{0}\} \mbox{ and } \Gh_{0}^{B}:= \{b_{0}\}.$$ Choose a decreasing sequence $(\delta_{i})_{i \in \N}$ of strictly positive numbers such that 
\begin{equation}
\label{34}
\sum_{i \in \N} \delta_{i} < \infty.
\end{equation} 
We may assume that $\delta_{0}=2$.

We construct the necessary data for \eqref{barbarphi-diagram} by induction. To begin, let 
\[
m_{-1} = \bar{m}_{0} = \bar{m}''_{-1} = n_{-1} = \bar{n}_{0}= \bar{n}''_{-1} := 0.
\]
Set 
\[
\bar{x}_{-1} = \be_{A \otimes Z_{\mathfrak{p},\mathfrak{q}}}
\]
and
\[
\bar{y}_{-1} = \be_{B \otimes Z_{\mathfrak{p},\mathfrak{q}}}.
\]
Define 
\[
\Gh_{-1}^{A} = \Gh_{-1}^{B} := \emptyset,
\]
\[
\Gh_{0}^{A}:= \{a_{0}\}
\]
and
\[
\Gh_{0}^{B}:= \{b_{0}\}.
\]

For the induction step, suppose that for some $i \in \N$ and $j=0,\ldots,i$, we have constructed $\Gh^{A}_{j}$, $\Gh^{B}_{j}$, $m_{j-1}$, $ \bar{m}_{j}$, $\bar{m}_{j-1}''$, $n_{j-1}$, $\bar{n}_{j}$, $\bar{n}_{j-1}''$, $\bar{x}_{j-1}$ and $\bar{y}_{j-1}$ with the following properties (which are easily verified for $i=0$):
\begin{enumerate} 
\item 
\[
\bar{m}''_{j-1} \ge m_{j-1}
\]
and
\[
\bar{x}_{j-1} \in A \otimes \gamma_{2\bar{m}''_{j-1}}(Z_{P_{2 \bar{m}''_{j-1}},Q_{2 \bar{m}''_{j-1}}}) \subset A \otimes Z_{\mathfrak{p},\mathfrak{q}}
\]
is a unitary satisfying
\[
\id_{A} \otimes \varrho_{m_{j-1}} =_{(\Gh_{j-1}^{A},3\delta_{j-1})} \ad(\bar{x}_{j-1}) \circ \varphi^{-1} \circ (\id_{B} \otimes \varrho_{n_{j-1}}) \circ \varphi \circ (\id_{A} \otimes \varrho_{m_{j-1}})
\]
\item
\begin{eqnarray*}
\Gh_{j}^{A}& := & \Gh_{j-1}^{A} \cup \{a_{j}\} \cup \varphi^{-1} \circ (\id_{B} \otimes \varrho_{n_{j-1}})(\Gh_{j-1}^{B})\\
& & \cup \bigcup_{k=0}^{j-1} (\ad(\bar{x}_{j-1}^{*})\circ (\id_{A} \otimes \varrho_{m_{j-1}}) \circ \ldots \circ \ad(\bar{x}_{k}^{*})\circ (\id_{A} \otimes \varrho_{m_{k}})(\Gh^{A}_{j-1}))\\
& & \subset A \otimes Z_{\mathfrak{p},\mathfrak{q}}
\end{eqnarray*}
\item $\bar{m}_{j}$ comes from Lemma~\ref{rho-m-rho-n-intertwining}, applied to $\Gh_{j}^{A}$ in place of $\Fh$ and $\delta_{j}$ in place of $\varepsilon$
\item
\[
m_{j}:= \max\{\bar{m}_{j}, \bar{m}''_{j-1},m_{j-1}+1\}
\]
\item
letting $m_{j}$ play the role of $m$, let $\bar{n}_{j} \in \N$ be the number (corresponding to $\bar{n}$) obtained from Lemma~\ref{rho-m-rho-n-intertwining}
\item
\[
\bar{n}''_{j-1} \ge n_{j-1} 
\]
and
\[
\bar{y}_{j-1} \in B \otimes \gamma_{2 \bar{n}''_{j-1}}(Z_{P_{2\bar{n}''_{j-1}},Q_{2\bar{n}''_{j-1}}})
\]
is a unitary satisfying
\[
\id_{B} \otimes \varrho_{n_{j-1}} = _{(\Gh^{B}_{j-1},3\delta_{j-1})} \ad(\bar{y}_{j-1}) \circ \varphi \circ  (\id_{A} \otimes \varrho_{m_{j}}) \circ \varphi^{-1} \circ (\id_{B} \otimes \varrho_{n_{j-1}})
\]
\item
\begin{eqnarray*}
\Gh_{j}^{B} & := & \Gh_{j-1}^{B} \cup \{b_{j}\} \cup \varphi \circ (\id_{A} \otimes \varrho_{m_{j}})(\Gh_{j}^{A})\\
& & \cup \bigcup_{k=0}^{j-1} (\ad(\bar{y}_{j-1}^{*})\circ (\id_{B} \otimes \varrho_{n_{j-1}}) \circ \ldots \circ \ad(\bar{y}_{k}^{*})\circ (\id_{B} \otimes \varrho_{n_{k}})(\Gh^{B}_{j-1})) \\
& & \subset B \otimes Z_{\mathfrak{p},\mathfrak{q}}.
\end{eqnarray*}
\end{enumerate}

We proceed to construct $\Gh^{A}_{i+1}$, $\Gh^{B}_{i+1}$, $m_{i+1}$, $ \bar{m}_{i+1}$, $\bar{m}_{i}''$, $n_{i}$, $\bar{n}_{i+1}$, $\bar{n}_{i}''$, $\bar{x}_{i}$ and $\bar{y}_{i}$.

Interchanging the roles of $A$ and $B$ (and replacing $\varphi$ by $\varphi^{-1}$), we may again apply Lemma~\ref{rho-m-rho-n-intertwining}  (this time with  $\Gh_{i}^{B}$ in place of $\Fh$ and $\delta_{i}$ in place of $\varepsilon$), to obtain a number $\bar{n}'_{i} \in \N$ (corresponding to the number $\bar{m}$ of \ref{rho-m-rho-n-intertwining}). Set 
\[
n_{i}:= \max\{\bar{n}_{i},\bar{n}'_{i},n_{i-1}+1\}.
\] 
Letting $n_{i}$ play the role of $m$, we obtain some $\bar{m}'_{i} \in \N$ (corresponding to the number $\bar{n}$ of \ref{rho-m-rho-n-intertwining}). We may clearly assume 
\begin{equation}
\label{11c}
\bar{m}'_{i} \ge m_{i}.
\end{equation}

By property (iii) above and by the assertion of Lemma~\ref{rho-m-rho-n-intertwining}, there is a unitary $x_{i} \in B \otimes \Zh$ such that 
\begin{equation}
\label{12c}
\varphi \circ (\id_{A} \otimes \varrho_{m_{i}}) =_{(\Gh_{i}^{A},\delta_{i})} \ad(x_{i}) \circ (\id_{B} \otimes \varrho_{n_{i}}) \circ \varphi \circ (\id_{A} \otimes \varrho_{m_{i}}).
\end{equation}

Since $\varphi^{-1}(x_{i}) \in A \otimes Z_{\mathfrak{p},\mathfrak{q}}$ is a unitary, there are  
\[
\bar{m}''_{i} \ge \bar{m}'_{i} \ge m_{i}
\]
and a unitary 
\begin{equation}
\label{18c}
\bar{x}_{i} \in A \otimes \gamma_{2\bar{m}''_{i}}(Z_{P_{2 \bar{m}''_{i}},Q_{2 \bar{m}''_{i}}}) \subset A \otimes Z_{\mathfrak{p},\mathfrak{q}}
\end{equation}
such that 
\begin{equation}
\label{13c}
\bar{x}_{i} =_{\delta_{i}} \varphi^{-1}(x_{i}),
\end{equation}
whence 
\begin{equation}
\label{27}
\id_{A} \otimes \varrho_{m_{i}} \stackrel{\eqref{12c},\eqref{13c}}{=_{(\Gh_{i}^{A},3\delta_{i})}} \ad(\bar{x}_{i}) \circ \varphi^{-1} \circ (\id_{B} \otimes \varrho_{n_{i}}) \circ \varphi \circ (\id_{A} \otimes \varrho_{m_{i}}).
\end{equation}

Set
\begin{eqnarray}
\label{14c}
\Gh_{i+1}^{A}& := &  \Gh_{i}^{A} \cup \{a_{i+1}\} \cup \varphi^{-1} \circ (\id_{B} \otimes \varrho_{n_{i}})(\Gh_{i}^{B}) \nonumber \\
& & \cup \bigcup_{k=0}^{i} (\ad(\bar{x}_{i}^{*})\circ (\id_{A} \otimes \varrho_{m_{i}}) \circ \ldots \circ (\ad(\bar{x}_{k}^{*})\circ (\id_{A} \otimes \varrho_{m_{k}})(\Gh^{A}_{i})).
\end{eqnarray}

Apply Lemma~\ref{rho-m-rho-n-intertwining}, this time with $\Gh_{i+1}^{A}$ in place of $\Fh$ and $\delta_{i+1}$ in place of $\varepsilon$, to obtain $\bar{m}_{i+1} \in \N$. Set
\begin{equation}
\label{15c}
m_{i+1}:= \max\{\bar{m}_{i+1}, \bar{m}''_{i},m_{i}+1\}.
\end{equation}
Letting $m_{i+1}$ play the role of $m$, let $\bar{n}_{i+1} \in \N$ be the number (corresponding to $\bar{n}$) obtained from Lemma~\ref{rho-m-rho-n-intertwining}.

Lemma~\ref{rho-m-rho-n-intertwining} (and our choice of $n_{i}$ and $m_{i+1}$ -- note that $m_{i+1} \ge \bar{m}'_{i}$) now yields a unitary 
\[
y_{i} \in A  \otimes Z_{\mathfrak{p},\mathfrak{q}}
\]
such that
\[
\varphi^{-1} \circ (\id_{B} \otimes \varrho_{n_{i}}) = _{(\Gh^{B}_{i},\delta_{i})} \ad(y_{i}) \circ (\id_{A} \otimes \varrho_{m_{i+1}}) \circ \varphi^{-1} \circ (\id_{B} \otimes \varrho_{n_{i}}).
\]
Similarly as before, there are 
\begin{equation}
\label{17c}
\bar{n}''_{i} \ge \bar{n}'_{i} \ge n_{i}
\end{equation}
and a unitary 
\begin{equation}
\label{19c}
\bar{y}_{i} \in B \otimes \gamma_{2\bar{n}''_{i}}(Z_{P_{2\bar{n}''_{i}},Q_{2\bar{n}''_{i}}})
\end{equation}
such that 
\[
\bar{y}_{i} =_{\delta_{i}} \varphi(y_{i}),
\]
whence
\begin{equation}
\label{30}
\id_{B} \otimes \varrho_{n_{i}} = _{(\Gh^{B}_{i},3\delta_{i})} \ad(\bar{y}_{i}) \circ \varphi \circ  (\id_{A} \otimes \varrho_{m_{i+1}}) \circ \varphi^{-1} \circ (\id_{B} \otimes \varrho_{n_{i}}).
\end{equation}
Set
\begin{eqnarray}
\label{16c}
\Gh_{i+1}^{B}
& := & \Gh_{i}^{B} \cup \{b_{i+1}\}  \cup  \varphi \circ(\id_{A} \otimes \varrho_{m_{i+1}})(\Gh_{i+1}^{A}) \nonumber \\
& & \cup \bigcup_{k=0}^{i} (\ad(\bar{y}_{i}^{*})\circ (\id_{B} \otimes \varrho_{n_{i}}) \circ \ldots \circ \ad(\bar{y}_{k}^{*})\circ (\id_{B} \otimes \varrho_{n_{k}})(\Gh^{B}_{i})).
\end{eqnarray}

This completes the induction step. Properties (i)...(vii) above (with $i$ replaced by $i+1$) follow from \eqref{11c}, \eqref{27}, \eqref{14c}, \eqref{15c}, \eqref{17c}, \eqref{30} and \eqref{16c}. With these data, we obtain a diagram 
\[
\xymatrix @C=-0.7em @W=3.5em @R7ex{
\cdots \ar[rr]& & A {\otimes} Z_{\mathfrak{p},\mathfrak{q}} \ar[dr]_(.4){\alpha_{i}} \ar[rr]^{\zeta_{i}} && A {\otimes} Z_{\mathfrak{p},\mathfrak{q}} \ar[dr] \ar[rr] &&\cdots &   \ar[rr]& & \bar{A} \ar@<-0.3em>[d]_{\bar{\varphi}} \\
&\cdots  \ar[rr] && B {\otimes} Z_{\mathfrak{p},\mathfrak{q}} \ar[ur]^(.4){\beta_{i}} \ar[rr]_{\xi_{i}} && B  {\otimes} Z_{\mathfrak{p},\mathfrak{q}} \ar[ur] \ar[rr] &&    \cdots \ar[rr] & & \bar{B}   \ar@<-0.3em>[u]_{\bar{\varphi}^{-1}}
}
\]
where
\begin{equation}
\label{28}
\zeta_{i}= \ad(\bar{x}_{i}^{*}) \circ (\id_{A} \otimes \varrho_{m_{i}})
\end{equation}
and
\begin{equation}
\label{29}
\xi_{i} = \ad(\bar{y}_{i}^{*})\circ (\id_{B} \otimes \varrho_{n_{i}}),
\end{equation}
which is an asymptotic intertwining in the sense of \cite{BlaKir:limits}. More precisely, the rows are generalized inductive systems of c.p.c.\ maps with limit maps $\zeta_{i,\infty}$ and $\xi_{i,\infty}$, respectively and   limits $\bar{A}$ and $\bar{B}$ (which are $C^{*}$-algebras). The maps $\bar{\varphi}$ and $\bar{\varphi}^{-1}$ are mutually inverse $*$-isomorphisms, given by 
\begin{eqnarray}
\label{barphi-limit}
\bar{\varphi} \circ \zeta_{i,\infty} & = & \lim_{j \to \infty} \xi_{j,\infty} \circ \alpha_{j} \circ \zeta_{j-1} \circ \ldots \circ \zeta_{i} 
\end{eqnarray}
and
\begin{eqnarray}
\label{barphiinverse-limit}
\bar{\varphi}^{-1} \circ \xi_{i,\infty} & = & \lim_{j \to \infty} \zeta_{j+1,\infty} \circ \beta_{j} \circ \xi_{j-1} \circ \ldots \circ \xi_{i} 
\end{eqnarray}
for any $i \in \N$. To check existence of these limits, consider $a \in \Gh_{j}^{A}$ and note that
\begin{eqnarray}
\label{33}
\lefteqn{
\|\alpha_{j+1} \zeta_{j}(a) - \xi_{j} \alpha_{j}(a)\|} \nonumber\\
& \le & \|\alpha_{j+1} \zeta_{j}(a) - \alpha_{j+1} \beta_{j} \alpha_{j}(a)\| + \|\alpha_{j+1} \beta_{j} \alpha_{j}(a) - \xi_{j} \alpha_{j}(a)\| \nonumber\\
& \le & 2 \cdot 3 \delta_{j},
\end{eqnarray}
where we have used \eqref{31}, \eqref{32}, \eqref{28}, \eqref{29}, \eqref{27} and \eqref{30} for the second estimate. Using \eqref{33} and \eqref{34}, induction now shows that the sequence
\[
\xi_{j,\infty} \circ \alpha_{j} \circ \zeta_{j-1} \circ \ldots \circ \zeta_{i} (d)
\]
is Cauchy (hence convergent in $B \otimes \Zh$) for any $d \in \Gh_{i}^{A}$. But then, the limit yields a linear map 
\[
\bigcup_{i \in \N} \mathrm{span } \{\zeta_{i,\infty}(\Gh_{i}^{A})\} \to B \otimes \Zh.
\]
This map is clearly contractive, hence extends to all of $\bar{A}$;  we denote the resulting map by $\bar{\varphi}$. A similar reasoning works for $\bar{\varphi}^{-1}$, and it is straightforward to check that these maps indeed are mutually inverse $*$-homomorphisms.   We have now shown  \eqref{barphi-limit} and \eqref{barphiinverse-limit}.

Next, we  define elements $u_{i} \in A \otimes Z_{\mathfrak{p},\mathfrak{q}}$ and $z_{i} \in B \otimes Z_{\mathfrak{p},\mathfrak{q}}$ for $i \in \N$ recursively by
\begin{equation}
\label{20c}
u_{0}:= \be_{A \otimes Z_{\mathfrak{p},\mathfrak{q}}} , \; u_{i}:= \bar{x}_{i-1}^{*} (\id_{A} \otimes \varrho_{m_{i-1}})(u_{i-1})
\end{equation}
and
\begin{equation}
\label{21c}
z_{0}:= \be_{B \otimes Z_{\mathfrak{p},\mathfrak{q}}} , \; z_{i}:= \bar{y}^{*}_{i-1} (\id_{B} \otimes \varrho_{n_{i-1}})(z_{i-1}).
\end{equation}
From \eqref{18c}, \eqref{19c} and \eqref{8c} above we see that $u_{i-1}$ and $z_{i-1}$ live in the multiplicative domains of the u.c.p.\ maps $\id_{A} \otimes \varrho_{m_{i-1}}$ and $\id_{B} \otimes \varrho_{n_{i-1}}$, respectively, whence $u_{i}$ and $z_{i}$ again are unitaries. 

We then obtain a diagram
\[
\xymatrix @C=-0.7em @W=3.5em @R7ex{
\cdots \ar[rr]& & A {\otimes} Z_{\mathfrak{p},\mathfrak{q}} \ar[d]_{\ad(u_{i})} \ar[rr]^{\id_{A} \otimes \varrho_{m_{i}}} && A {\otimes} Z_{\mathfrak{p},\mathfrak{q}} \ar[d]_{\ad(u_{i+1})} \ar[rr] &&\cdots &   \ar[rr]& & A \otimes \Zh \ar[d]_{\cong} \ar@/^3ex/[ddd]^{\bar{\bar{\varphi}}}\\
\cdots \ar[rr]& & A {\otimes} Z_{\mathfrak{p},\mathfrak{q}} \ar[dr]_(.4){\alpha_{i}} \ar[rr]^{\zeta_{i}} && A {\otimes} Z_{\mathfrak{p},\mathfrak{q}} \ar[dr] \ar[rr] &&\cdots &   \ar[rr]& & \bar{A} \ar[d]_{\bar{\varphi}} \\
&\cdots  \ar[rr] && B {\otimes} Z_{\mathfrak{p},\mathfrak{q}} \ar[ur]^(.4){\beta_{i}} \ar[rr]_{\xi_{i}} \ar[d]_{\ad(z_{i}^{*})} && B  {\otimes} Z_{\mathfrak{p},\mathfrak{q}} \ar[ur] \ar[rr] \ar[d]_{\ad(z_{i+1}^{*})} &&    \cdots \ar[rr] & & \bar{B} \ar[d]_{\cong} \\
&\cdots  \ar[rr] && B {\otimes} Z_{\mathfrak{p},\mathfrak{q}}  \ar[rr]_{\id_{B} \otimes \varrho_{n_{i}}} && B  {\otimes} Z_{\mathfrak{p},\mathfrak{q}}  \ar[rr] &&    \cdots \ar[rr] & & B \otimes \Zh
}
\]
in which the upper and the lower thirds are unitary intertwinings. Commutativity of the rectangles follows from \eqref{28}, \eqref{29}, \eqref{20c} and \eqref{21c}. Using \eqref{barphi-limit}, we see that the isomorphism $\bar{\bar{\varphi}}$ is given by 
\begin{eqnarray*}
\lefteqn{\bar{\bar{\varphi}} \circ \ldots \circ (\id_{A} \otimes \varrho_{m_{i+1}}) \circ (\id_{A} \otimes \varrho_{m_{i}}) }\\
& = & \lim_{j \to \infty} \ldots \circ (\id_{B} \otimes   \varrho_{n_{j}})  \ad(z_{j}^{*})  \alpha_{j}  \ad(u_{j})  (\id_{A} \otimes \varrho_{m_{j-1}})  \circ \ldots \circ (\id_{A} \otimes \varrho_{m_{i}})
\end{eqnarray*}
for each $i \in \N$; we obtain a similar statement for $\bar{\bar{\varphi}}^{-1}$. Setting
\[
w^{B}_{i}:=\ldots \circ (\id_{B} \otimes \varrho_{n_{i+1}}) \circ (\id_{B} \otimes \varrho_{n_{i}})(z_{i}^{*} \alpha_{i}(u_{i}))
\]
and
\[
w^{A}_{i}:=  \ldots \circ (\id_{A} \otimes \varrho_{m_{i+2}})\circ (\id_{A} \otimes \varrho_{m_{i+1}})  (u_{i+1}^{*} \beta_{i}(z_{i}))
\]
for $i \in \N$, it is now straightforward to check property (i) of the proposition.

To verify (ii), note that (i) in particular implies 
\[
\bar{\bar{\varphi}} \circ (\id_{A} \otimes \be_{\Zh}) = \lim_{i \to \infty} \ad(w^{B}_{i}) \circ ( \id_{B} \otimes (\ldots \circ \varrho_{n_{i+1}} \circ \varrho_{n_{i}})) \circ \varphi \circ (\id_{A} \otimes \be_{Z_{\mathfrak{p},\mathfrak{q}}}).
\]
By Definition~\ref{def-p-q-exponential-system}(iv), we have a sequence $(v_{i})_{i \in \N} \subset \Zh$ of unitaries such that 
\[
\ad(v_{i}) \circ (\ldots \circ \varrho_{n_{i+1}} \circ \varrho_{n_{i}}) \stackrel{i \to \infty}{\longrightarrow} \bar{\sigma}_{\mathfrak{p},\mathfrak{q}}
\]
pointwise. But then, we obtain
\begin{eqnarray*}
\lefteqn{\bar{\bar{\varphi}} \circ (\id_{A} \otimes \be_{\Zh})}\\
& = & \lim_{i \to \infty} \ad(w^{B}_{i}) \circ \ad(\be_{B} \otimes v_{i}^{*}) \circ (\id_{B} \otimes \bar{\sigma}_{\mathfrak{p},\mathfrak{q}}) \circ \varphi \circ (\id_{A} \otimes \be_{Z_{\mathfrak{p},\mathfrak{q}}}),
\end{eqnarray*}
the desired approximate unitary equivalence.
\end{nproof}
\en

\bn
\label{suspended-iso-invariant}
\begin{nprop}
Let $\mathfrak{p}$ and $\mathfrak{q}$ be supernatural numbers of infinite type which are relatively prime. Suppose $A$ and $B$ are separable, unital, $\Zh$-stable $C^{*}$-algebras and let 
\[
\varphi: A  \otimes Z_{\mathfrak{p},\mathfrak{q}} \to B  \otimes Z_{\mathfrak{p},\mathfrak{q}}
\]
be a unitarily suspended $\Ch([0,1])$-isomorphism. Then, there is an isomorphism 
\[
\tilde{\varphi}: A  \to B  \otimes \Zh
\]
such that
\begin{equation}
\label{23c}
\tilde{\varphi} \au (\id_{B} \otimes \bar{\sigma}_{\mathfrak{p},\mathfrak{q}}) \circ \varphi \circ (\id_{A} \otimes \be_{Z_{\mathfrak{p},\mathfrak{q}}}),
\end{equation}
where $\bar{\sigma}_{\mathfrak{p},\mathfrak{q}}$ is the standard embedding from Proposition~\ref{Z-p-q-embedding}.
\end{nprop}

\begin{nproof}
Since $A$ and $B$ are both $\Zh$-stable, we may as well replace them by $A \otimes \Zh$ and $B \otimes \Zh$, respectively. 

Apply Proposition~\ref{barbarphi-intertwining} (with $A\otimes \Zh$ and $B\otimes \Zh$ in place of $A$ and $B$, respectively) to obtain an isomorphism 
\[
\bar{\bar{\varphi}}: A \otimes \Zh \otimes \Zh \to B \otimes \Zh \otimes \Zh.  
\]
Since $\Zh$ is strongly self-absorbing, there is an isomorphism $\theta: \Zh \to \Zh \otimes \Zh$ such that 
\begin{equation}
\label{22c}
\theta \au \id_{\Zh} \otimes \be_{\Zh}.
\end{equation} 
Set
\[
\tilde{\varphi}:= \bar{\bar{\varphi}} \circ (\id_{A} \otimes \theta): A \otimes \Zh \to B \otimes \Zh \otimes \Zh.
\]
Then, $\tilde{\varphi}$ is an isomorphism and we have  
\[
\tilde{\varphi} \au \bar{\bar{\varphi}} \circ (\id_{A \otimes \Zh} \otimes \be_{\Zh});
\]
the assertion follows from Proposition~\ref{barbarphi-intertwining}(ii).
\end{nproof}
\en

\bn
\label{p-phi-barphi}
\begin{nprop}
Let $\mathfrak{p}$ and $\mathfrak{q}$ be supernatural numbers of infinite type which are relatively prime. Suppose $A$ and $B$ are separable, unital, $\Zh$-stable $C^{*}$-algebras and let 
\[
\varphi: A  \otimes Z_{\mathfrak{p},\mathfrak{q}} \to B  \otimes Z_{\mathfrak{p},\mathfrak{q}}
\]
be a unitarily suspended $\Ch([0,1])$-isomorphism. Let 
\[
\Lambda= (\lambda,\Delta): \Inv (A) \to \Inv(B)
\]
be an isomorphism of invariants. Suppose that
\begin{equation}
\label{phi-T-Delta}
\varphi^{T}(\tau \otimes \bar{\tau}_{Z_{\mathfrak{p},\mathfrak{q}}}) = \Delta(\tau) \otimes \bar{\tau}_{Z_{\mathfrak{p},\mathfrak{q}}} \; \forall \, \tau \in T(B)
\end{equation}
(where $\bar{\tau}_{Z_{\mathfrak{p},\mathfrak{q}}}$ is the canonical tracial state on $Z_{\mathfrak{p},\mathfrak{q}}$ induced by the Lebesgue measure on $[0,1]$) and that the diagram
\begin{equation}
\label{A-B-lambda-diagram}
\xymatrix{
K_{*}(A) \otimes K_{*}(Z_{\mathfrak{p},\mathfrak{q}}) \ar[r]^(.55){\alpha} \ar[d]_{\lambda \otimes \id} & K_{*}(A \otimes Z_{\mathfrak{p},\mathfrak{q}}) \ar[d]^{\varphi_{*}} \\
K_{*}(B) \otimes K_{*}(Z_{\mathfrak{p},\mathfrak{q}}) \ar[r]^(.55){\alpha} & K_{*}(B \otimes Z_{\mathfrak{p},\mathfrak{q}})
}
\end{equation}
commutes (where the maps $\alpha$ come from the K{\"u}nneth theorem, cf.\ \cite[Section~23]{Bla:newk-theory}). Then, there is an isomorphism 
\[
\bar{\varphi}: A  \to B 
\]
with
\[
\Inv(\bar{\varphi}) = \Lambda.
\]
\end{nprop}

\begin{nproof}
Again, since $A$ and $B$ are both $\Zh$-stable, we may as well replace them by $A \otimes \Zh$ and $B \otimes \Zh$, respectively. Note that \eqref{phi-T-Delta} then translates to
\begin{equation}
\label{25c}
\varphi^{T}(\tau \otimes \bar{\tau}_{\Zh} \otimes \bar{\tau}_{Z_{\mathfrak{p},\mathfrak{q}}}) = \Delta(\tau \otimes \bar{\tau}_{\Zh}) \otimes \bar{\tau}_{Z_{\mathfrak{p},\mathfrak{q}}} \; \forall \, \tau \in T(B).
\end{equation}
Let 
\[
\tilde{\varphi}: A \otimes \Zh \to B \otimes \Zh \otimes \Zh
\]
be the isomorphism from Proposition~\ref{suspended-iso-invariant} and set 
\[
\bar{\varphi} := (\id_{B} \otimes \theta^{-1}) \circ \tilde{\varphi},
\]
where $\theta: \Zh \stackrel{\cong}{\longrightarrow} \Zh \otimes \Zh$ is as in the proof of \ref{suspended-iso-invariant}. Consider the diagram
\begin{equation}
\xymatrix@R=3em{
K_{*}(A{\otimes}\Zh) \ar[r]^{=} \ar[d]^{(\id_{A {\otimes}\Zh})_{*} {\otimes} [\be_{Z_{\mathfrak{p},\mathfrak{q}}}]} & K_{*}(A{\otimes}\Zh) \ar[r]^{=} \ar[d]^{(\id_{A {\otimes} \Zh} {\otimes} \be_{Z_{\mathfrak{p},\mathfrak{q}}})_{*}} &K_{*}(A{\otimes}\Zh) \ar[r]^{=} \ar[d]^{(\id_{A {\otimes} \Zh} {\otimes} \be_{\Zh})_{*}}&K_{*}(A{\otimes}\Zh) \ar[d]^{=} \\
K_{*}(A {\otimes} \Zh) {\otimes} K_{*}(Z_{\mathfrak{p},\mathfrak{q}}) \ar[r]^(.55){\alpha} & K_{*}(A {\otimes} \Zh {\otimes} Z_{\mathfrak{p},\mathfrak{q}}) \ar[r]^{(\id {\otimes} \bar{\sigma}_{\mathfrak{p},\mathfrak{q}})_{*}} & K_{*}(A {\otimes} \Zh {\otimes} \Zh) \ar[r]^(.55){(\id{\otimes}\theta)^{-1}_{*}} & K_{*}(A {\otimes} \Zh).
}
\label{A-otimes-Z-diagram}
\end{equation}
Identifying the upper left copy of $K_{*}(A \otimes \Zh)$ with $K_{*}(A \otimes \Zh) \otimes K_{*}(\C)$, we see from the K{\"u}nneth theorem that the first rectangle commutes. The second rectangle commutes, since  it even does so at the level of algebras.

The third rectangle commutes, since $\theta^{-1} \circ (\id_{\Zh} \otimes \be_{\Zh}) \au \id_{\Zh}$ by \eqref{22c}. Thus, we see that \eqref{A-otimes-Z-diagram} commutes; the same holds for the respective diagram with $A$ replaced by $B$.   The first vertical map is an isomorphism by Proposition~\ref{KK-C-Z-p-q}, whence we obtain a commuting diagram of isomorphisms
\begin{equation}
\label{A-B-otimes-Z-diagram}
\xymatrix@R=3em{
K_{*}(A{\otimes}\Zh) \ar[r]_(.38){\id {\otimes} [\be]}  \ar[d]^{\lambda} \ar@/^4ex/[rrr]^{=} & K_{*}(A{\otimes}\Zh) {\otimes} K_{*}(Z_{\mathfrak{p},\mathfrak{q}}) \ar[r]_(.55){\alpha} \ar[d]^{\lambda \otimes \id} & K_{*}(A{\otimes}\Zh {\otimes}Z_{\mathfrak{p},\mathfrak{q}}) \ar[r]_(.55){(\id {\otimes} \theta^{-1} \bar{\sigma}_{\mathfrak{p},\mathfrak{q}})_{*}} & K_{*}(A{\otimes}\Zh) \ar[d]^{\lambda}\\
K_{*}(B{\otimes}\Zh) \ar[r]^(.38){\id {\otimes} [\be]} \ar@/_4ex/[rrr]_{=} & K_{*}(B{\otimes}\Zh) {\otimes} K_{*}(Z_{\mathfrak{p},\mathfrak{q}}) \ar[r]^(.55){\alpha} & K_{*}(B{\otimes}\Zh {\otimes}Z_{\mathfrak{p},\mathfrak{q}}) \ar[r]^(.55){(\id {\otimes} \theta^{-1} \bar{\sigma}_{\mathfrak{p},\mathfrak{q}})_{*}} & K_{*}(B{\otimes}\Zh).
} 
\end{equation}
On the other hand, by \eqref{23c} we have 
\begin{equation}
\label{24c}
\bar{\varphi} \au (\id_{B} \otimes \theta^{-1}) \circ (\id_{B \otimes \Zh} \otimes \bar{\sigma}_{\mathfrak{p},\mathfrak{q}}) \circ \varphi \circ (\id_{A \otimes \Zh} \otimes \be_{\mathfrak{p},\mathfrak{q}}),
\end{equation}
and from commutativity of the first rectangle of \eqref{A-otimes-Z-diagram} we see that the diagram
\begin{equation}
\label{incomplete-diagram}
\xymatrix@R=3em{
K_{*}(A{\otimes}\Zh) \ar[r]_(.38){\id {\otimes} [\be]}   \ar@/^4ex/[rrr]^{=} & K_{*}(A{\otimes}\Zh) {\otimes} K_{*}(Z_{\mathfrak{p},\mathfrak{q}}) \ar[r]_(.55){\alpha} & K_{*}(A{\otimes}\Zh {\otimes}Z_{\mathfrak{p},\mathfrak{q}})   \ar[d]^{\varphi_{*}} & K_{*}(A{\otimes}\Zh) \ar[d]^{\bar{\varphi}_{*}}\\
 &  & K_{*}(B{\otimes}\Zh {\otimes}Z_{\mathfrak{p},\mathfrak{q}}) \ar[r]^(.55){(\id {\otimes} \theta^{-1} \bar{\sigma}_{\mathfrak{p},\mathfrak{q}})_{*}} & K_{*}(B{\otimes}\Zh).
}
\end{equation}
commutes. Putting together \eqref{incomplete-diagram}, \eqref{A-B-otimes-Z-diagram} and \eqref{A-B-lambda-diagram}, we obtain a commuting diagram
\[
\xymatrix{K_{*}(A \otimes \Zh) \ar[r]^{=} \ar[d]_{\lambda} & K_{*}(A \otimes \Zh) \ar[d]^{\bar{\varphi}_{*}} \\
K_{*}(B \otimes \Zh) \ar[r]^{=}& K_{*}(B \otimes \Zh).
}
\]
Since $\lambda$ and $\bar{\varphi}_{*}$ are isomorphisms, we in fact have
\[
\lambda = \bar{\varphi}_{*}.
\]
Note that \eqref{incomplete-diagram}, \eqref{A-B-otimes-Z-diagram} and \eqref{A-B-lambda-diagram} a priori only show that $\lambda$ and $\bar{\varphi}_{*}$ agree as isomorphisms of abelian groups. However, since they are both isomorphisms of ordered abelian groups with order unit, they also have to agree as such. 

Finally, for any $\tau \in T(B)$ we have
\begin{eqnarray*}
\bar{\varphi}^{T}(\tau \otimes \bar{\tau}_{\Zh}) & \stackrel{\eqref{24c}}{=} & ((\id_{B} \otimes \theta^{-1}) \circ (\id_{B \otimes \Zh} \otimes \bar{\sigma}) \circ \varphi \circ (\id_{A \otimes \Zh} \otimes \be_{Z_{\mathfrak{p},\mathfrak{q}}}))^{T}(\tau \otimes \bar{\tau}_{\Zh}) \\
& = & (\id_{A} \otimes \be_{Z_{\mathfrak{p},\mathfrak{q}}})^{T} \circ \varphi^{T} \circ (\id_{B \otimes \Zh} \otimes \bar{\sigma}_{\mathfrak{p},\mathfrak{q}})^{T} \circ (\id_{B} \otimes \theta^{-1})^{T}(\tau \otimes \bar{\tau}_{\Zh}) \\
& \stackrel{\ref{Z-p-q-embedding}}{=} & (\id_{A \otimes \Zh} \otimes \be_{Z_{\mathfrak{p},\mathfrak{q}}})^{T} \circ \varphi^{T} (\tau \otimes \bar{\tau}_{\Zh} \otimes \bar{\tau}_{\Zh_{\mathfrak{p},\mathfrak{q}}}) \\
&\stackrel{\eqref{25c}}{=} & (\id_{A \otimes \Zh} \otimes \be_{Z_{\mathfrak{p},\mathfrak{q}}})^{T} (\Delta(\tau \otimes \bar{\tau}_{\Zh}) \otimes \bar{\tau}_{\Z_{\mathfrak{p},\mathfrak{q}}}) \\
& = & \Delta(\tau \otimes \bar{\tau}_{\Zh}),
\end{eqnarray*}
whence
\[
\Delta = \bar{\varphi}^{T} .
\]
Therefore, 
\[
\Lambda =  (\lambda,\Delta)  = (\bar{\varphi}_{*},\bar{\varphi}^{T}) = \Inv (\bar{\varphi}),
\]
as desired.
\end{nproof}
\en

\bn
The results of this section motivate the following definition; it  was also  inspired by a suggestion from H.\ Lin of how to repair a mistake in an initial  version of Corollary~\ref{fewtraces-lfdr-cor} below. 

\label{lift-Lambda-along-Z}
\begin{ndefn}
Let $A$ and $B$ be separable, simple, unital and nuclear $C^{*}$-algebras and let $\mathfrak{p}$ and  $\mathfrak{q}$ be supernatural numbers which are relatively prime. We say an isomorphism $\Lambda: \Inv(A) \to \Inv(B)$ can be lifted along $Z_{\mathfrak{p},\mathfrak{q}}$, if there is a unitarily suspended $\Ch([0,1])$-isomorphism $\varphi: A \otimes Z_{\mathfrak{p},\mathfrak{q}} \to B \otimes Z_{\mathfrak{p},\mathfrak{q}}$ such that (with the notation of Definition~\ref{d-unitarily-suspended} and Proposition~\ref{lambda-lambda-U} below) 
\[
\Inv(\sigma_{\mathfrak{p}}) = \Lambda_{M_{\mathfrak{p}}} \mbox{ and } \Inv(\varrho_{\mathfrak{q}}) = \Lambda_{M_{\mathfrak{q}}}.
\]
\end{ndefn}
\en

\section{$K$-theory of tensor products with UHF algebras}

Our proof of Theorem~\ref{main-result} relies on the analysis of tensor products with generalized dimension drop intervals. These are $\Ch([0,1])$-algebras with each fibre being a tensor product with a UHF algebra. Below we collect some easy observations about the $K$-theory of such algebras.

\bn
\label{lambda-lambda-U}
\begin{nprop}
Let $A$ and $B$ be separable, simple, unital $C^{*}$-algebras; suppose 
\[
\Lambda=(\lambda,\Delta): \Inv(A) \to \Inv(B)
\]
is an isomorphism of invariants. For an infinite dimensional UHF algebra $\Uh$, let 
\[
\lambda_{\Uh}:= \alpha_{B,\Uh} \circ (\lambda \otimes \id) \circ \alpha_{A,\Uh}^{-1} 
\]
be the isomorphism (of abelian groups)  that makes 
\[
\xymatrix{
K_{*}(A) \otimes K_{*}(\Uh) \ar[r]^{\alpha_{A,\Uh}} \ar[d]_{\lambda \otimes \id} & K_{*}(A \otimes \Uh) \ar[d]^{\lambda_{\Uh}}\\
K_{*}(B) \otimes K_{*}(\Uh) \ar[r]^{\alpha_{B,\Uh}} & K_{*}(B \otimes \Uh) 
}
\]
commute (where $\alpha_{A,\Uh}$ and $\alpha_{B,\Uh}$ are the isomorphisms obtained from the K\"unneth theorem). Let
\[
\Delta_{\Uh}: T(B \otimes \Uh) \to T(A \otimes \Uh)
\]
be given by 
\begin{equation}
\label{26c}
\Delta_{\Uh}(\tau \otimes \bar{\tau}_{\Uh}) := \Delta(\tau) \otimes \bar{\tau}_{\Uh},
\end{equation}
where $\bar{\tau}_{\Uh} \in T(\Uh)$ is the unique tracial state. Then, 
\[
\Lambda_{\Uh}:= (\lambda_{\Uh},\Delta_{\Uh})
\]
is an isomorphism of invariants. (In particular, $\lambda_{\Uh}$ is an isomorphism of ordered $K$-groups.) 
\end{nprop}

\begin{nproof}
The K{\"u}nneth maps $\alpha_{A,\Uh}$ and $\alpha_{B,\Uh}$ send elements of the form $[p]\otimes [q]$ to elements of the form $[p \otimes q]$ (up to identifying $\Kh \otimes \Kh$ with $\Kh$).

This in particular implies  that 
\[
(\tau \otimes \bar{\tau}_{\Uh})_{*} \circ \alpha_{B,\Uh} = \tau_{*} \otimes (\bar{\tau}_{\Uh})_{*}
\]
for any trace $\tau \in T(B)$. The respective statement holds for $A$, so 
\[
\Delta(\tau)_{*} = (\Delta(\tau)_{*} \otimes (\bar{\tau}_{\Uh})_{*}) \circ \alpha_{A,\Uh}^{-1}.
\]
Now if $x \in K_{*}(A \otimes \Uh)$ and $\tau \in T(B)$, we obtain
\begin{eqnarray}
\label{23}
(\tau \otimes \bar{\tau}_{\Uh})_{*} \circ \lambda_{\Uh}(x) & = & (\tau \otimes \bar{\tau}_{\Uh})_{*} \circ \alpha_{B,\Uh} \circ (\lambda \otimes \id) \circ \alpha_{A,\Uh}^{-1} (x) \nonumber \\
& = & (\tau_{*} \circ \lambda \otimes (\bar{\tau}_{\Uh})_{*}) \circ \alpha_{A,\Uh}^{-1} (x) \nonumber \\
& = & (\Delta(\tau)_{*} \otimes (\bar{\tau}_{\Uh})_{*}) \circ \alpha_{A,\Uh}^{-1} (x) \nonumber \\
& = & (\Delta(\tau) \otimes \bar{\tau}_{\Uh})_{*} (x) \nonumber \\
& \stackrel{\eqref{26c}}{=} & (\Delta_{\Uh}(\tau \otimes \bar{\tau}_{\Uh}))_{*}(x).
\end{eqnarray}
By \cite[Theorem~5(a)]{JiaSu:Z}, $A \otimes \Uh$ and $B \otimes \Uh$ are $\Zh$-stable, so the order structure on $K_{0}$ is determined by traces, cf.\ \cite{Ror:Z-absorbing}; by simplicity this just means that $0 \neq x \in K_{0}$ is positive if and only if $\tau_{*}(x) >0$ for all traces. Since $\Delta_{\Uh}$ is a bijection, \eqref{23} then implies that $\lambda_{\Uh}$  maps positive elements to positive elements, hence is an isomorphism of ordered abelian groups; the same holds for $\lambda^{-1}_{\Uh}$. At the same time, we see from \eqref{23} that
\[
\xymatrix{
T(A \otimes \Uh) \ar[r]^{r_{A \otimes \Uh}} & S(K_{0}(A \otimes \Uh))  \\
T(B\otimes \Uh) \ar[r]^{r_{B \otimes \Uh}} \ar[u]^{\Delta_{\Uh}} & S(K_{0}(B \otimes \Uh)) \ar[u]_{\lambda_{\Uh}^{S}}
}
\]
commutes, so $\Lambda_{\Uh}$ indeed is an isomorphism of invariants in the sense of Definition~\ref{morphisms-of-invariants}.
\end{nproof}
\en

\bn
\label{p-plus-q-injectivity}
The proof of the following observation was pointed out to us by M.\ R{\o}rdam.

\begin{nprop}
Let $\mathfrak{p}$ and $\mathfrak{q}$ be supernatural numbers which are relatively prime. Suppose $G$ is an  abelian group. Then, the homomorphism 
\[
\beta: G \to G \otimes K_{0}(M_{\mathfrak{p}}) \oplus G \otimes K_{0}(M_{\mathfrak{q}})
\]
given by 
\[
\beta = \id_{G} \otimes [\be_{M_{\mathfrak{p}}}] \oplus \id_{G} \otimes [\be_{M_{\mathfrak{q}}}]
\]
is injective. In particular, if $A$ is a separable $C^{*}$-algebra, the composition of maps
\[
\xymatrix{
K_{*}(A) \ar[r]^{\id {\otimes} \be}_{\cong} &K_{*}(A {\otimes} Z_{\mathfrak{p},\mathfrak{q}})  \ar[r]^{\alpha^{-1}(\ev_{0} \oplus \ev_{1})_{*}}  & K_{*}(A) {\otimes} (K_{*}(M_{\mathfrak{p}})  {\oplus} K_{*}(M_{\mathfrak{q}})),
}
\]
where 
\[
\alpha: K_{*}(A) {\otimes} (K_{*}(M_{\mathfrak{p}})  {\oplus} K_{*}(M_{\mathfrak{q}})) \to K_{*}(A {\otimes} M_{\mathfrak{p}} \oplus A {\otimes} M_{\mathfrak{q}})
\]
is the isomorphism obtained from the K{\"u}nneth theorem, is injective.
\end{nprop}

\begin{nproof}
Writing $\beta$ as  $\beta_{\mathfrak{p}} \oplus \beta_{\mathfrak{q}} $ in the obvious way, it will suffice to show that
\begin{equation}
\label{25}
\ker \beta_{\mathfrak{p}} \subset \{g \in G \mid \exists p \in \N :  p|\mathfrak{p}, \, pg = 0\}  
\end{equation}
and 
\begin{equation}
\label{26}
\ker \beta_{\mathfrak{q}} \subset \{g \in G \mid \exists q \in \N :  q|\mathfrak{q}, \, qg = 0\}.
\end{equation}
For then, it is clear that 
\[
\ker \beta_{\mathfrak{p}} \cap \ker \beta_{\mathfrak{q}} = \{0\}
\]
since $\mathfrak{p}$ and $\mathfrak{q}$ are relatively prime. 

To prove \eqref{25}, note that  $G \otimes K_{0}(M_{\mathfrak{p}})$ can be written as an inductive limit
\[
G \otimes \Z \stackrel{\mu_{1,2}}{\longrightarrow} G \otimes \Z \stackrel{\mu_{2,3}}{\longrightarrow} \ldots  \longrightarrow G \otimes K_{0}(M_{\mathfrak{p}}),
\]
where each connecting map $\mu_{i,i+1}$ is just multiplication by some integer $p_{i}$ which divides $\mathfrak{p}$. In this picture, the map $\beta_{\mathfrak{p}}$ is given by $\mu_{1,\infty} \circ (\id_{G} \otimes 1)$. 
Now if $g \in G$ is such that $pg \neq 0$ for all $p \in \N$ dividing $\mathfrak{p}$, each $\mu_{1,i}(g \otimes 1)$ is nonzero, whence $\beta_{\mathfrak{p}}(g) = \mu_{1,\infty}(g \otimes 1) \neq 0$ and $g \notin \ker \beta_{\mathfrak{p}}$.   So, we have proven \eqref{25}; \eqref{26} is verified in the same way. 
\end{nproof}
\en

\section{Tracial rank zero and an asymptotic unitary equivalence}

The purpose of this section is to recall the notion of tracial rank zero and to show how existing classification results  can be used to obtain a classification result for tensor products with generalized dimension drop algebras, provided all the fibre algebras have tracial rank zero.

\bn
\label{d-tr0}
In \cite{Lin:TAF1}, Lin introduced the notion of tracially AF algebras -- or, equivalently, $C^{*}$-algebras of tracial rank zero. In the simple and unital case, the definition reads as follows:
 
\begin{ndefn}
A simple, unital $C^{*}$-algebra $A$ is said to be tracially AF, if, for any finite subset $\Fh \subset A$, $\varepsilon > 0$ and $0 \neq a \in A_{+}$, there is a finite-dimensional $C^{*}$-subalgebra $B \subset A$ with the following properties:
\begin{enumerate}
\item[(i)] $\|b \be_{B} - \be_{B} b\|<\varepsilon \; \forall \, b \in \Fh$ 
\item[(ii)] $\dist(\be_{B}b\be_{B}, B)< \varepsilon \; \forall \, b \in \Fh$ 
\item[(iii)] $\be_{A}-\be_{B}$ is Murray--von Neumann equivalent to a projection in the hereditary subalgebra $\overline{aAa}$ of $A$.
\end{enumerate}
\end{ndefn}
\en

\bn
\label{sigma-rho-asu}
The following result from \cite{Lin:existence} (which is based on \cite{Lin:asu-auto}, \cite{Lin:appr-hom-hom} and \cite{KisKum:Ext}) will play a crucial role for the proof of \ref{fewtraces-lfdr-cor}. 

\begin{ntheorem}
Let $A$ and $B$ be separable, simple, unital and nuclear $C^{*}$-algebras and let $\mathfrak{p}$ and  $\mathfrak{q}$ be supernatural numbers of infinite type which are relatively prime and which satisfy $M_{\mathfrak{p}} \otimes M_{\mathfrak{q}} \cong \Qh$. Suppose that $A \otimes M_{\mathfrak{p}}$, $A \otimes M_{\mathfrak{q}}$, $B \otimes M_{\mathfrak{p}}$ and $B \otimes M_{\mathfrak{q}}$ have tracial rank zero, and that either $K_{*}(A \otimes M_{\mathfrak{p}})$ is torsion-free, or $K_{*}(A)$ contains its torsion part as a direct summand.  Let 
\[
\sigma_{\mathfrak{p}}: A \otimes M_{\mathfrak{p}} \to B \otimes M_{\mathfrak{p}}
\]
and 
\[
\varrho_{\mathfrak{q}}: A \otimes M_{\mathfrak{q}} \to B \otimes M_{\mathfrak{q}}
\]
be two unital isomorphisms. For the obvious induced maps 
\[
\sigma,\varrho: A \otimes M_{\mathfrak{p}} \otimes M_{\mathfrak{q}} \to B \otimes M_{\mathfrak{p}} \otimes M_{\mathfrak{q}}
\]
assume 
\[
[\sigma] = [\varrho] \in KK(A \otimes \Qh,B \otimes \Qh).
\]
Then, there is an isomorphism 
\[
\sigma'_{\mathfrak{p}}: A \otimes M_{\mathfrak{p}} \to B \otimes M_{\mathfrak{p}}
\]
such that
\[
\Inv(\sigma'_{\mathfrak{p}}) = \Inv(\sigma_{\mathfrak{p}})
\]
and such that the induced isomorphism
\[
\sigma': A \otimes M_{\mathfrak{p}} \otimes M_{\mathfrak{q}} \to B \otimes M_{\mathfrak{p}} \otimes M_{\mathfrak{q}}
\]
is asymptotically unitarily equivalent to $\varrho$. 
\end{ntheorem}
\en

\bn
\label{lifting-along-Z-p-q-cor}
\begin{ncor}
Let $A$ and $B$ be  separable, simple, unital $C^{*}$-algebras with locally finite decomposition rank,  satisfying the UCT; suppose that, for $A$ and $B$, projections separate traces. Furthermore, suppose that  the torsion part of $K_{*}(A)$ misses at least one  prime order, or that $K_{*}(A)$ contains its torsion part as a direct summand. 
 
Then, for any isomorphism $\Lambda: \Inv(A) \to \Inv(B)$, there are supernatural numbers of infinite type $\mathfrak{p}$ and $\mathfrak{q}$ which are relatively prime, which satisfy $M_{\mathfrak{p}} \otimes M_{\mathfrak{q}} \cong \Qh$ and such that $\Lambda$ can be lifted along $Z_{\mathfrak{p},\mathfrak{q}}$.
\end{ncor}

\begin{nproof}
Since $K_{*}(A)$ misses a prime order, it is easy to find supernatural numbers of infinite type $\mathfrak{p}$ and $\mathfrak{q}$ which are relatively prime, which satisfy $M_{\mathfrak{p}} \otimes M_{\mathfrak{q}} \cong \Qh$ and such that $K_{*}(A \otimes M_{\mathfrak{p}})$ is torsion-free.  Note that  $A \otimes M_{\mathfrak{p}}$ is approximately divisible and hence has real rank zero by \cite{BlaKumRor:apprdiv} since projections separate traces.   Now by \cite[Theorem~2.1]{Winter:lfdr}, $A \otimes M_{\mathfrak{p}}$ has tracial rank zero; a similar reasoning works for $B \otimes M_{\mathfrak{p}}$ etc. One may now apply Theorem~\ref{sigma-rho-asu} to obtain the desired unitarily suspended isomorphism which lifts $\Lambda$ along $Z_{\mathfrak{p},\mathfrak{q}}$. In the case where $K_{*}(A)$ contains its torsion part as a direct summand, Theorem~\ref{sigma-rho-asu} may be applied in the same way.
\end{nproof}
\en

\section{Localizing at the Jiang--Su algebra}

\bn
\label{main-result}
\begin{ntheorem}
Let $\Ah$ be a class of separable, simple, unital and nuclear $C^{*}$-algebras such that, for any $A$ and $B$ in $\Ah$ and  any isomorphism of invariants $\Lambda: \Inv(A) \to \Inv(B)$, there are relatively prime supernatural numbers of infinite type $\mathfrak{p}$ and $\mathfrak{q}$ for which $\Lambda$ can be lifted along $Z_{\mathfrak{p},\mathfrak{q}}$. 
 
Then, $\Ah$ satisfies  the Elliott conjecture localized at $\Zh$.
\end{ntheorem}

\begin{nproof}
Let $A$ and  $B$ be in $\Ah^{\Zh}$ (i.e., $A$ and $B$ are both in $\Ah$ and are $\Zh$-stable); suppose $$\Lambda =(\lambda,\Delta): \Inv(A) \to \Inv(B)$$ is an isomorphism of the Elliott invariants.  We have to establish an isomorphism between $A$ and $B$ which induces $\Lambda$. 

By hypothesis, there are supernatural numbers of infinite type $\mathfrak{p}$ and $\mathfrak{q}$ which are relatively prime and such that $\Lambda$ can be lifted along $Z_{\mathfrak{p},\mathfrak{q}}$.

In view of Proposition~\ref{p-phi-barphi}, it will suffice to show that for the unitarily suspended  isomorphism 
\[
\varphi: A  \otimes Z_{\mathfrak{p},\mathfrak{q}} \to B  \otimes Z_{\mathfrak{p},\mathfrak{q}}
\]
lifting $\Lambda$ we have
\begin{equation}
\label{phi-Delta}
\varphi^{T}(\tau  \otimes \bar{\tau}_{Z_{\mathfrak{p},\mathfrak{q}}}) = \Delta(\tau ) \otimes \bar{\tau}_{Z_{\mathfrak{p},\mathfrak{q}}} \; \forall \, \tau \in T(B)
\end{equation}
and that the diagram
\begin{equation}
\label{lambda-phi-diagram}
\xymatrix{
K_{*}(A) \otimes K_{*}(Z_{\mathfrak{p},\mathfrak{q}}) \ar[r]^{\alpha} \ar[d]_{\lambda \otimes \id} & K_{*}(A \otimes Z_{\mathfrak{p},\mathfrak{q}}) \ar[d]^{\varphi_{*}} \\
K_{*}(B) \otimes K_{*}(Z_{\mathfrak{p},\mathfrak{q}}) \ar[r]^{\alpha} & K_{*}(B \otimes Z_{\mathfrak{p},\mathfrak{q}}) .
}
\end{equation}
commutes. For then, by Proposition~\ref{p-phi-barphi} there will be an isomorphism $\bar{\varphi}:A \to B$ such that $\Inv(\bar{\varphi})= \Lambda$.

 Let 
\[
\Lambda_{\mathfrak{p}}:= (\lambda_{\mathfrak{p}},\Delta_{\mathfrak{p}}): \Inv(A \otimes M_{\mathfrak{p}}) \to \Inv(B \otimes M_{\mathfrak{p}}),
\]
\[
\Lambda_{\mathfrak{q}}:= (\lambda_{\mathfrak{q}},\Delta_{\mathfrak{q}}): \Inv(A \otimes M_{\mathfrak{q}}) \to \Inv(B \otimes M_{\mathfrak{q}}),
\]
and
\[
\Lambda_{\Qh}:= (\lambda_{\Qh},\Delta_{\Qh}): \Inv(A \otimes M_{\Qh}) \to \Inv(B \otimes M_{\Qh}),
\]
be the induced isomorphisms of invariants  from Proposition~\ref{lambda-lambda-U}.

By hypothesis,  there are  isomorphisms  
\[
\sigma_{\mathfrak{p}}: A \otimes M_{\mathfrak{p}} \to B \otimes M_{\mathfrak{p}}
\]
and 
\[
\varrho_{\mathfrak{q}}: A \otimes M_{\mathfrak{q}} \to B \otimes M_{\mathfrak{q}}
\]
such that 
\begin{equation}
\label{sigma-p-rho-q}
\Inv(\sigma_{\mathfrak{p}}) = \Lambda_{\mathfrak{p}} \mbox{ and } \Inv(\varrho_{\mathfrak{q}}) = \Lambda_{\mathfrak{q}}.
\end{equation}
From this we obtain isomorphisms 
\[
\sigma,\varrho : A \otimes M_{\mathfrak{p}} \otimes M_{\mathfrak{q}} \to B \otimes M_{\mathfrak{p}}\otimes M_{\mathfrak{q}},
\]
where $\sigma$ is just $\sigma_{\mathfrak{p}} \otimes \id_{M_{\mathfrak{q}}}$ and $\varrho$ coincides with     $\varrho_{\mathfrak{q}}^{[1,3]} \otimes \id_{M_{\mathfrak{p}}}^{[2]}$ (i.e., up to identifying $M_{\mathfrak{q}} \otimes M_{\mathfrak{p}}$ with $M_{\mathfrak{p}} \otimes M_{\mathfrak{q}}$ in the obvious way).  

Using the K\"unneth theorem, it follows that 
\begin{equation}
\label{sigma-rho-lambdaQ}
\sigma_{*} = \varrho_{*} = \lambda_{\Qh} : K_{*}(A \otimes \Qh) \to K_{*}(B \otimes \Qh)
\end{equation}
(the K\"unneth theorem a priori only yields equality of isomorphisms between abelian groups -- but since $\lambda_{\Qh}$ is an isomorphism of the \emph{ordered} $K$-groups, so are $\sigma_{*}$ and $\varrho_{*}$). It is trivial to check that
\[
\sigma^{T} = \varrho^{T} = \Delta_{\Qh};
\]
therefore, we have 
\begin{equation}
\label{Inv-Lambda}
\Inv(\sigma) = \Lambda_{\Qh} = \Inv(\varrho).
\end{equation}

For $t \in [0,1]$, let 
\[
\varphi_{t}:(A \otimes Z_{\mathfrak{p},\mathfrak{q}})_{t}  \to (B \otimes Z_{\mathfrak{p},\mathfrak{q}})_{t}
\]
be the fibre map of $\varphi$. It is clear from \eqref{77} and \eqref{Inv-Lambda}, that  
\[
\varphi_{t}^{T}(\tau \otimes \bar{\tau}_{t}) = \Delta(\tau) \otimes \bar{\tau}_{t}
\]
for any $\tau \in T(B)$, where $\bar{\tau}_{t}$ denotes the unique tracial state on the UHF algebra $(Z_{\mathfrak{p},\mathfrak{q}})_{t}$. But then,  \eqref{phi-Delta} holds by Proposition~\ref{p-Z-p-q-tau-preserving}.

To check that diagram \eqref{lambda-phi-diagram} commutes, note that
\begin{equation*}
\xymatrix{
K_{*}(A)  \ar[r]^{(\id_{A})_{*} \otimes [\be_{Z_{\mathfrak{p},\mathfrak{q}}}]} \ar[d]_{\lambda} & K_{*}(A) \otimes K_{*}(Z_{\mathfrak{p},\mathfrak{q}}) \ar[d]^{\lambda \otimes \id} \\
K_{*}(B)  \ar[r]^{(\id_{B})_{*} \otimes [\be_{Z_{\mathfrak{p},\mathfrak{q}}}]} & K_{*}(B) \otimes K_{*}(Z_{\mathfrak{p},\mathfrak{q}}) 
}
\end{equation*}
commutes by the K{\"u}nneth theorem and by Proposition~\ref{KK-C-Z-p-q}, whence it suffices to check that 
\begin{equation}
\label{22}
\xymatrix{
K_{*}(A)  \ar[r]^{(\id_{A} \otimes \be_{Z_{\mathfrak{p},\mathfrak{q}}})_{*}} \ar[d]_{\lambda} & K_{*}(A \otimes Z_{\mathfrak{p},\mathfrak{q}}) \ar[d]^{\varphi_{*}} \\
K_{*}(B)  \ar[r]^{(\id_{B} \otimes \be_{Z_{\mathfrak{p},\mathfrak{q}}})_{*}} & K_{*}(B \otimes Z_{\mathfrak{p},\mathfrak{q}}) 
}
\end{equation}
commutes. To this end, note that 
\begin{equation}
\label{20}
\xymatrix{
K_{*}(A {\otimes} Z_{\mathfrak{p},\mathfrak{q}}) \ar[d]^{\varphi_{*}} \ar[r]^{(\ev_{0} \oplus \ev_{1})_{*}} & K_{*}(A {\otimes} M_{\mathfrak{p}} \oplus A {\otimes} M_{\mathfrak{q}}) \ar[d]^{(\sigma_{\mathfrak{p}} \oplus \varrho_{\mathfrak{q}})_{*}} \ar[r]^{\alpha^{-1}} & K_{*}(A) {\otimes} (K_{*}(M_{\mathfrak{p}})  {\oplus} K_{*}(M_{\mathfrak{q}})) \ar[d]^{\lambda {\otimes} (\id \oplus  \id)} \\
K_{*}(B {\otimes} Z_{\mathfrak{p},\mathfrak{q}})   \ar[r]^{(\ev_{0} \oplus \ev_{1})_{*}} & K_{*}(B {\otimes} M_{\mathfrak{p}} \oplus A {\otimes} M_{\mathfrak{q}})  \ar[r]^{\alpha^{-1}} & K_{*}(B) {\otimes} ( K_{*}(M_{\mathfrak{p}})  {\oplus} K_{*}(M_{\mathfrak{q}}))
}
\end{equation}
(using \eqref{sigma-p-rho-q}, \eqref{77} and the K{\"u}nneth theorem) and 
\begin{equation}
\label{21}
\xymatrix{
K_{*}(A)  \ar[r]^{(\id_{A} {\otimes} \be_{Z_{\mathfrak{p},\mathfrak{q}}})_{*}} \ar[d]_{\lambda} & K_{*}(A {\otimes} Z_{\mathfrak{p},\mathfrak{q}})  \ar[r]^{\alpha^{-1}(\ev_{0} \oplus \ev_{1})_{*}} & K_{*}(A) {\otimes} (K_{*}(M_{\mathfrak{p}})  \oplus K_{*}(M_{\mathfrak{q}})) \ar[d]^{\lambda {\otimes} (\id \oplus  \id)} \\
K_{*}(B)  \ar[r]^{(\id_{B} {\otimes} \be_{Z_{\mathfrak{p},\mathfrak{q}}})_{*}} & K_{*}(B {\otimes} Z_{\mathfrak{p},\mathfrak{q}}) \ar[r]^{\alpha^{-1}(\ev_{0} \oplus \ev_{1})_{*}} & K_{*}(B) {\otimes} ( K_{*}(M_{\mathfrak{p}})  \oplus K_{*}(M_{\mathfrak{q}}))
}
\end{equation}
 (using naturality of the K{\"u}nneth theorem) both commute. Now by Proposition~\ref{p-plus-q-injectivity}, the horizontal lines of \eqref{20} and \eqref{21} are injective. Combining \eqref{20} and \eqref{21} we see that \eqref{22} indeed commutes, from which in turn follows  that \eqref{lambda-phi-diagram} commutes.

Since we have already verified \eqref{phi-Delta}, this proves the Elliott conjecture for $\Ah$, localized at $\Zh$.
\end{nproof}
\en

\section{Applications}

\bn
\label{fewtraces-lfdr-cor}
As a consequence of Theorem~\ref{main-result} and the results of \cite{Lin:existence}, we obtain the following far-reaching generalization of \cite[Corollary~2.3]{Winter:lfdr}. We are indebted to H.\ Lin for pointing out a mistake in an earlier version of the corollary below, and for suggesting a way to repair the proof. 

\begin{ncor}
Let $\Ah$ be the class of separable, simple, unital $C^{*}$-algebras with locally finite decomposition rank, and satisfying the UCT. Suppose that, for each $A \in \Ah$, projections separate traces. Furthermore, suppose that for each $A \in \Ah$ either the torsion part of $K_{*}(A)$ misses at least one  prime order, or that $K_{*}(A)$ contains its torsion part as a direct summand.  Then, $\Ah$ satisfies the   Elliott conjecture localized at $\Zh$.
\end{ncor}

\begin{nproof}
By Corollary~\ref{lifting-along-Z-p-q-cor}, $\Ah$ satisfies the hypothesis of Theorem~\ref{main-result}.
\end{nproof}
\en

\bn
\label{monotracial-cor}
In the monotracial case, we obtain a version of \ref{fewtraces-lfdr-cor} for limits of type I $C^{*}$-algebras:

\begin{ncor}
The Elliott conjecture, localized at $\Zh$, holds for the class of all separable, simple, unital inductive limits of type I $C^{*}$-algebras with unique tracial state and which satisfy any of the $K$-theory restrictions of \ref{fewtraces-lfdr-cor}. 
\end{ncor}

\begin{nproof}
Let $\Ah$ denote the class of the corollary and let $A$ be in $\Ah$.  For a supernatural number $\mathfrak{p}$ of infinite type, $A \otimes M_{\mathfrak{p}}$ has real rank zero, stable rank one and weakly unperforated $K_0$-group by \cite{BlaKumRor:apprdiv}. Therefore, $A \otimes M_{\mathfrak{p}}$ has tracial rank zero by \cite[Corollary~5.16]{Lin:traces}.  It is well-known that limits of type I $C^{*}$-algebras satisfy the UCT, so we may apply our classification result in the form of \cite[Theorem~3.5]{Lin:existence} to conclude that $\Ah^{\Zh}$ satisfies the Elliott conjecture. 
\end{nproof}
\en

\bn
\label{ASH-cor}
We note the following special case of \ref{fewtraces-lfdr-cor} explicitly:

\begin{ncor}
The Elliott conjecture, localized at $\Zh$, holds for the class of all separable, simple, unital ASH algebras with   finitely generated $K$-theory for which projections separate traces.  
\end{ncor}

\begin{nproof}
By \cite[Corollary~2.2]{NgWinter:subhom},  any separable, unital ASH algebra $A$ has locally finite decomposition rank.  In fact, by \cite[Corollary~2.1]{NgWinter:subhom} (cf.\ also \cite[Theorem~1.6]{Winter:subhomdr}), $A$ may be written as a limit of recursive subhomogeneous algebras $B_{n}$ of finite topological dimension. By \cite[Theorem~2.16]{Phillips:recsub}, any such $B_{n}$ is an iterated pullback of algebras of the form $\Ch(X) \otimes M_{r}$. But then, $A$ clearly satisfies the UCT. Since the  $K$-theory of $A$ is finitely generated, so is its torsion part. We have thus checked that the class of the corollary is contained in the class $\Ah$ of Corollary~\ref{fewtraces-lfdr-cor}.
\end{nproof}
\en

\bn
\label{uniquely-ergodic-cor}
As a particularly striking application of the preceding corollaries, we may use results of Q.\ Lin and Phillips (\cite{LinPhi:mindifflimits}) to make progress on the classification problem for minimal dynamical systems:

\begin{ncor}
Up to $\Zh$-stability, $C^{*}$-algebras of uniquely ergodic, minimal, smooth dynamical systems are classified by their ordered $K$-theory.
\end{ncor}
\en

\bn
\label{Z-ASH-characterization}
Employing a recent result of  Dadarlat and  Toms, we obtain a very satisfactory characterization of $\Zh$ among limits of type I  $C^{*}$-algebras:

\begin{ncor}
The Jiang--Su algebra $\Zh$ is the uniquely determined strongly self-absorbing and projectionless limit of type I $C^{*}$-algebras.
\end{ncor}

\begin{nproof}
Let $\Dh$ be a strongly self-absorbing limit  of type I $C^{*}$-algebras. It is well-known that $\Dh$ satisfies the UCT. Now by the results of \cite[Section~5]{TomsWinter:ssa}, $K_{0}(\Dh)=\Z=K_{0}(\Zh)$ and $K_{1}(\Dh)=0=K_{1}(\Zh)$. By a recent result of Dadarlat and Toms (cf.\ \cite{DadToms:Z-ASH-uniqueness}),  $\Dh$ absorbs the Jiang--Su algebra. It has unique trace by \cite[Theorem~1.7]{TomsWinter:ssa}. Since $\Dh$ is $\Zh$-stable, the order structure on $K_{0}(\Dh)$ is determined by its trace -- which clearly implies that the Elliott invariants of $\Dh$ and $\Zh$ agree. Hence, $\Dh \cong \Zh$ by Corollary~\ref{monotracial-cor}.
\end{nproof}
\en

\bn
\label{ASH-range}
Using results of Elliott and Thomsen, we can describe the classes of Corollaries~\ref{fewtraces-lfdr-cor} and \ref{monotracial-cor} more precisely, at least up to $\Zh$-stability:

\begin{ncor}
Let $\Ah$ denote the union of the classes of Corollaries~\ref{fewtraces-lfdr-cor} and \ref{monotracial-cor}. Then, the class $\Ah^{\Zh}$ consists entirely of ASH algebras of topological dimension at most 3. 

Moreover,  $A \in \Ah^{\Zh}$ has real rank zero if and only if  $K_{0}(A)$ has Riesz interpolation and $K_{0}(A)/\mathrm{Tor}(K_{0}(A)) \neq \Z$; in this case, $A$ is AH of topological dimension at most 3. 
\end{ncor}

\begin{nproof}
By \cite{GongJiangSu:Z}, simple $\Zh$-stable $C^{*}$-algebras have weakly unperforated $K$-theory. By \cite[Theorem~3.4.4 and Example~3.4.8]{Ror:encyc} in connection with \cite[Example~1.10]{Winter:subhomdr},  for any $A \in \Ah$ there is a simple, unital ASH algebra $B$ of topological dimension at most 2 which has the same invariant as $A \otimes \Zh$. Corollary~\ref{fewtraces-lfdr-cor}  (or \ref{monotracial-cor}, respectively) show that $A \otimes \Zh \cong B \otimes \Zh$. 

$B \otimes \Zh$ clearly is ASH. Write $B$ as an inductive limit of subhomogeneous algebras $B_{i}$ each of which has topological dimension at most 2. Write $\Zh$ as limit of prime dimension drop intervals $Z_{p_{i},q_{i}}$. We may choose the $p_{i}$ and $q_{i}$ to grow arbitrarily fast, and  such that $p_{i} \ll q_{i} \ll p_{i+1}$ for all $i$. It is then straightforward to arrange that for any $i,k \in \N$ the space $\Prim_{k} (B_{i} \otimes Z_{p_{i},q_{i}})$   is either of the form $\Prim_{k'} B_{i}$ or $(\Prim_{k'} B_{i}) \times (0,1)$ (where $k'$ is either $k/p_{i}$, or $k/q_{i}$ or $k/(p_{i}q_{i}))$, respectively), whence $\dim (\Prim_{k} (B_{i} \otimes Z_{p_{i},q_{i}})) \le \dim (\Prim_{k'} B_{i}) + 1 \le 3$, cf.\ \cite{Winter:subhomdr}. This shows that the topological dimension of $B \otimes \Zh$ is at most 3.

If $A \in \Ah^{\Zh}$ has real rank zero,  it has tracial rank zero by \cite{Winter:lfdr}. By \cite{Lin:TAF1}, $K_{0}(A)$ has Riesz interpolation. 

If $A \in \Ah^{\Zh}$ has Riesz interpolation and $K_{0}(A)/\mathrm{Tor}(K_{0}(A)) \neq \Z$, by \cite{EllGong:rrzeroII} there is a simple AH algebra $B$ of topological dimension at most 3 which has real rank zero and the same invariant as $A$. By \cite{EllGongLi:apprdiv}, $B$ is approximately divisible, hence $\Zh$-stable by \cite{TomsWinter:Zash}. Our classification result yields $A \cong B$. (Note that, in this case, $A$ itself has real rank zero.) 
\end{nproof}
\en

\bn
\label{classified-examples}
\begin{nexamples}
The class of Corollary~\ref{fewtraces-lfdr-cor} in particular contains  the irrational rotation algebras, the (infinite dimensional) UHF algebras and the Bunce--Deddens algebras, which are all known to be $\Zh$-stable and to have real rank zero (cf.~\cite{Ror:encyc}). In fact, it covers all presently known separable, unital, simple, nuclear and $\Zh$-stable $C^{*}$-algebras which are finite and have real rank zero (and finitely generated $K$-theory). It  also contains the (projectionless) Jiang--Su algebra, and the (projectionless) crossed products of the form $\Ch(S^{3}) \rtimes_{\alpha} \Z$ considered in \cite[Section~5, Example~4]{Con:Thom}, with $\alpha$ a uniquely ergodic  minimal diffeomorphism (cf.\ \cite{LinPhi:mindifflimits}).
\end{nexamples}
\en

\bn
\begin{nremarks}
(i) It is remarkable that Corollary~\ref{fewtraces-lfdr-cor} does not depend on an inductive limit structure in any way; locally finite decomposition rank is a \emph{local} condition. 

(ii) In the monotracial case of Corollary~\ref{monotracial-cor}, the result does not in any way depend on the existence (or non-existence) of projections. As we have pointed out above, the class contains  projectionless $C^{*}$-algebras as well as examples with real rank zero.

(iii) The $K$-theory restriction of \ref{fewtraces-lfdr-cor}  does  not seem to be necessary. The technical reason why we could not remove it at the current stage is that \cite{Lin:existence} relies on \cite{KisKum:Ext}, which only covers AT-algebras (and the $K$-theory of these is always torsion free). 

(iv) The formulation \cite[Theorem~3.5]{Lin:existence} of our classification result shows that, in the situation of \ref{fewtraces-lfdr-cor}, the tracial state space does not have to be included in the invariant. (The reason is that $T(A)$ agrees with $T(A \otimes \Uh)$ for any UHF algebra $\Uh$, and that in the real rank zero case the tracial state space does not have to be included in the invariant.) This on the one hand means that our result does not go beyond the real rank zero case too far. On the other hand it lends credence to a point we made earlier, namely that, outside the real rank zero case, our result is perpendicular to Elliott--Gong--Li classification: What makes the class of Elliott--Gong--Li more general than real rank zero AH algebras is essentially the additional freedom coming from the map $r_{A}: T(A) \to S(K_{0}(A))$ being not injective.    

(v) In  Corollary~\ref{ASH-range}, we could only show that our classification result ranges over ASH algebras of topological dimension at most 3. However, we conjecture that topological dimension at most 2 will suffice. To prove this, it would be enough to show that the algebras of \cite[Theorem~3.4.4]{Ror:encyc} are $\Zh$-stable.    (In fact, we even conjecture that all  simple ASH algebras with finite topological dimension are $\Zh$-stable.) 
\end{nremarks}
\en

\section{Outlook}

\bn
\label{outlook}
Although Corollary~\ref{fewtraces-lfdr-cor} and its descendants look more prominent at the current stage, we nontheless like to think of Theorem~\ref{main-result}  as the main result of these notes, for two reasons. 

For once, \ref{main-result} implements an abstract method to pass from classification results with many projections to results with only few projections. The former are in general much easier to obtain, simply because  in the presence of sufficiently many projections, it is easier to read off information about the algebra from the Elliott invariant. 

The other reason is that Corollary~\ref{fewtraces-lfdr-cor} remains incomplete because of the condition on projections separating traces, and the $K$-theory restriction. Both  conditions should certainly not be necessary, and Theorem~\ref{main-result} together with the proof of \ref{fewtraces-lfdr-cor} indeed provide a clear strategy of how to remove them. To deal with the trace condition, in the proof of \ref{fewtraces-lfdr-cor}  replace `tracial rank zero' by `tracial rank one' (or by `TAI'), and verify the hypotheses of \ref{main-result}. First, one would have to generalize the results of \cite{Winter:lfdr} to tensor products with UHF algebras (but without asking for real rank zero).  One would then have to carry over the methods of \cite{Lin:existence} to the case where the target algebra is TAI (as opposed to TAF). To remove the $K$-theory restrictions, one would have to generalize \cite{KisKum:Ext} to more general classes of AH algebras. Although these are  technically hard problems, it is particularly promising that it now seems possible to use the actual classification result  of Elliott, Gong and Li rather than just their methods (which are bound to become much more inaccessible in the ASH case), to classify ASH algebras. 

If the above strategy works, the result will be the classification of separable,  simple, unital  $C^{*}$-algebras with locally finite decomposition rank and UCT up to $\Zh$-stability. At the current stage, this would  be the best such result in the simple, unital and stably finite case that is within sight. Replacing the hypothesis of locally finite decomposition rank by just stable finiteness would certainly be most desirable -- but it seems unlikely to the author that this question is within the scope of the methods that have been developed so far.   In fact, this may be closely  related to the UCT question itself.
\en

\bibliographystyle{amsplain}

\providecommand{\bysame}{\leavevmode\hbox to3em{\hrulefill}\thinspace}
\providecommand{\MR}{\relax\ifhmode\unskip\space\fi MR }
\providecommand{\MRhref}[2]{%
  \href{http://www.ams.org/mathscinet-getitem?mr=#1}{#2}
}
\providecommand{\href}[2]{#2}

\end{document}